\documentclass[10pt]{amsart}

\usepackage[colorlinks]{hyperref}
\usepackage{color,graphicx,shortvrb}
\usepackage[latin 1]{inputenc}

\usepackage[active]{srcltx} 

\usepackage{enumerate}
\usepackage{amssymb}
\usepackage{tikz-cd}

\usepackage [ all ]{xy}

\newtheorem{theorem}{Theorem}[section]

\newtheorem{corollary}[theorem]{Corollary}

\newtheorem{definition}[theorem]{Definition}
\newtheorem{lemma}[theorem]{Lemma}

\newtheorem{proposition}[theorem]{Proposition}
\newtheorem{remark}[theorem]{Remark}

\def\J#1#2#3{ \left\{ #1,#2,#3 \right\} }
\def\RR{{\mathbb{R}}}
\def\NN{{\mathbb{N}}}
\def\11{\textbf{$1$}}
\def\11b#1{\mathbf{1}_{_{#1}}}

\DeclareGraphicsExtensions{.jpg,.pdf,.png,.eps}

\begin{document}

\title[Surjective isometries between unitary sets of JB$^*$-algebras]{Surjective isometries between unitary sets of unital JB$^*$-algebras}

\author[M. Cueto-Avellaneda]{Mar{\'i}a Cueto-Avellaneda}

\address[M. Cueto-Avellaneda]{Departamento de An{\'a}lisis Matem{\'a}tico, Facultad de
Ciencias, Universidad de Granada, 18071 Granada, Spain.}
\curraddr{School of Mathematics, Statistics and Actuarial Science, University of Kent, Canterbury, Kent CT2 7NX, UK}
\email{mcueto@ugr.es, emecueto@gmail.com}

\author[Y. Enami]{Yuta Enami}
\address[Y. Enami]{Graduate School of Science and Technology, Niigata University, Niigata 950-2181, Japan}
\email{enami@m.sc.niigata-u.ac.jp}

\author[D. Hirota]{Daisuke Hirota}

\address[D. Hirota]{Graduate School of Science and Technology, Niigata University, Niigata 950-2181, Japan}
\email{f20j008d@mail.cc.niigata-u.ac.jp}

\author[T. Miura]{Takeshi Miura}

\address[T. Miura]{Department of Mathematics, Faculty of Science, Niigata University, Niigata 950-2181, Japan}
\email{miura@math.sc.niigata-u.ac.jp}

\author[A.M. Peralta]{Antonio M. Peralta}

\address[A.M. Peralta]{Departamento de An{\'a}lisis Matem{\'a}tico, Facultad de
Ciencias, Universidad de Granada, 18071 Granada, Spain.}
\email{aperalta@ugr.es}




\date{}

\begin{abstract} This paper is, in a first stage, devoted to establish a topological--algebraic characterization of the principal component, $\mathcal{U}^0 (M)$,  of the set of unitary elements, $\mathcal{U} (M)$, in a unital JB$^*$-algebra $M$. We arrive to the conclusion that, as in the case of unital C$^*$-algebras, $$\begin{aligned}\mathcal{U}^0(M) &=  M^{-1}_{\textbf{1}}\cap\mathcal{U} (M)
=\left\lbrace  U_{e^{i h_n}}\cdots U_{e^{i h_1}}(\textbf{1}) \colon \begin{array}{c}
                                                                      n\in \mathbb{N}, \ h_j\in M_{sa}  \\
                                                                      \forall\ 1\leq j \leq n
                                                                    \end{array}
  \right\rbrace \\
&= \left\lbrace  u\in \mathcal{U} (M) : \hbox{ there exists } w\in \mathcal{U}^0(M) \hbox{ with } \|u-w\|<2 \right\rbrace
\end{aligned}$$ is analytically arcwise connected. Actually, $\mathcal{U}^0(M)$ is the smallest quadratic subset of $\mathcal{U} (M)$ containing the set $e^{i M_{sa}}$. Our second goal is to provide a complete description of the surjective isometries between the principal components of two unital JB$^*$-algebras $M$ and $N$. Contrary to the case of unital C$^*$-algebras, we shall deduce the existence of connected components in $\mathcal{U} (M)$ which are not isometric as metric spaces. We shall also establish necessary and sufficient conditions to guarantee that a surjective isometry $\Delta: \mathcal{U}(M)\to \mathcal{U} (N)$ admits an extension to a surjective linear isometry between $M$ and $N$, a conclusion which is not always true. Among the consequences it is proved that $M$ and $N$ are Jordan $^*$-isomorphic if, and only if, their principal components are isometric as metric spaces if, and only if, there exists a surjective isometry $\Delta: \mathcal{U}(M)\to \mathcal{U}(N)$ mapping the unit of $M$ to an element in $\mathcal{U}^0(N)$. These results provide an extension to the setting of unital JB$^*$-algebras of the results obtained by O. Hatori for unital C$^*$-algebras.
\end{abstract}

\maketitle
\thispagestyle{empty}

\section{Introduction}

Why are we so attracted by the unitary group of unital C$^*$-algebras? The celebrated Russo--Dye theorem, asserting that the convex hull of the unitary elements in a unital C$^*$-algebra is norm dense in its closed unit ball (cf. \cite{RuDye}), is probably one of the eldest tools leading our attention to the unitary group. It is perhaps unnecessary to recall that an element, $u$, in a unital C$^*$-algebra, $A$, lies in the subgroup, $\mathcal{U} (A),$ of all unitaries in $A$ if $u u^* = u^* u =1$. As we have already advanced, $\mathcal{U} (A)$ is a subgroup of $A$ containing the unit element. Applications of the Russo-Dye theorem and subsequent generalizations have been employed along the last fifty five years, recent usages on problems about extensions of isometries appear, for example, in \cite{MoriOza2020} and \cite{BeCuFerPe2018, CuePer19}.\smallskip

The unitary group of a unital C$^*$-algebra was employed as a complete invariant to classify von Neumann factors and certain unital C$^*$-algebras. Highlighting the group structure of the set of unitaries, S. Sakai established that if $A$ and $B$ are AW$^*$-factors and $\rho: \mathcal{U} (A)\to \mathcal{U} (B)$ is a uniformly continuous group isomorphism between their respective unitary groups, then there is a unique map $T$ from $A$ onto $B$ which is either a linear or conjugate linear $^*$-isomorphism and which agrees with $\rho$ on $\mathcal{U} (A)$ (see \cite{Sak55}). A similar conclusion was proved by H.A. Dye for W$^*$-factors in \cite{Dye55} with an independent argument which also gives results when no continuity assumption is made on the group isomorphism. More recently, Al-Rawashdeh, Booth and Giordano show that if the unitary groups of two simple unital AH-algebras of slow dimension growth and of real rank zero are isomorphic as abstract groups, then their $K_0$-ordered groups are isomorphic (cf. \cite{AlRawBooGior2012}).\smallskip

In a surprising turn, O. Hatori and L. Moln{\'a}r considered the set $\mathcal{U}(A)$, of unitary elements in a unital C$^*$-algebra $A$, as a metric space equipped with the metric given by the C$^*$-norm and forget about its group structure (see \cite{HatMol2014}). It should be noted that the metric space $(\mathcal{U}(A),\|\cdot\|)$ is almost never compact (cf. \cite[Exercise II.2.1]{Tak}). In the case of von Neumann algebras, the metric space given by the set of unitaries is a complete invariant since, as shown by Hatori and Moln{\'a}r, every surjective isometry between the unitary groups of two von Neumann algebras admits an extension to a surjective real linear isometry between these algebras (see \cite[Corollary 3]{HatMol2014}). Actually if $\Delta: \mathcal{U} (A) \to \mathcal{U} (B)$ is a surjective isometry between the sets of unitary elements of two unital C$^*$-algebras, then there is a central projection $p \in B$ and a Jordan $^*$-isomorphism $J : A\to B$ satisfying $$\Delta (e^{ix}) = \Delta(1) (p J(e^{ix}) + (1-p) J(e^{ix})^*),$$ for all $x\in A_{sa}$ (see \cite[Theorem 1]{HatMol2014}). In particular $A$ and $B$ are Jordan $^*$-isomorphic. Here, and henceforth, the self-adjoint part of a C$^*$-algebra $A$ will be denoted by $A_{sa}$.\smallskip

The framework of von Neumann algebras is a very favorable scenario since, among other things, the set of unitaries in a von Neumann algebra $W$ is precisely the set $\exp(i W_{sa})$, which is an analytically arcwise connected set  (cf. \cite[Theorem 5.2.5]{KR1}). However, it is well known that in the case of a general unital C$^*$-algebra $A$ the set $\mathcal{U}(A)$ is not always connected (cf. \cite[Exercise I.11.3]{Tak},  \cite[Exercises 4.6.6 and 4.6.7]{KR1}, \cite[pages 162-164]{HatMol2014}, \cite{Hatori14}). The connected component of $\mathcal{U} (A)$ containing the unit element is called the principal component, and it is denoted by $\mathcal{U}^0(A)$.\smallskip 

Hatori and Moln{\'a}r pointed out in \cite[Corollary 8]{HatMol2014} the existence of surjective isometries between the sets of unitaries in two unital commutative C$^*$-algebras which do not admit an extension to a surjective real linear isometry between these algebras. The problem of determining the structure of all surjective isometries between the sets of unitaries in two arbitrary unital C$^*$-algebras was carried out by O. Hatori in \cite{Hatori14}, where two main results are established: Suppose $A$ and $B$ are two unital C$^*$-algebras. \begin{enumerate}[$(a)$]
	\item Suppose $\Delta: \mathcal{U}^0 (A) \to \mathcal{U}^0 (B)$ is a mapping. Then $\Delta$ is a surjective isometry if and only if there exist a central projection $p \in B$ and a Jordan $^*$-isomorphism $J : A\to B$ satisfying $$\Delta (e^{ix}) = \Delta(1) \Big(p J(a) + (1-p) J(a)\Big),$$ for all $a\in \mathcal{U}^0 (A)$. In particular $A$ and $B$ are Jordan $^*$-isomorphic. The principal components  $\mathcal{U}^0 (A)$ and $\mathcal{U}^0 (B)$ can be somehow replaced by any two connected components (see \cite[Theorem 3.1 and Corollary 3.5]{Hatori14}).
\item Each surjective isometry $\Delta: \mathcal{U} (A) \to \mathcal{U} (B)$ can be written as a direct sum of a family of surjective real linear isometries between $A$ and $B$ restricted to the corresponding connected components, and the extensibility of $\Delta$ to a surjective real linear isometry is only possible under additional conditions (\cite[Theorem 4.1 and Corollary 5.1]{Hatori14}). 
\end{enumerate}

From a strict mathematical point of view, C$^*$-algebras are contained in the strictly wider class of JB$^*$-algebras, a class  of Jordan-Banach algebras determined by approprite version of the Gelfand-Naimark axiom. The class of JB$^*$-algebras is  nowadays a consolidated object of study, whose origins go back to the works of P. Jordan, J. von Neumann, E. Wigner, I. Kaplansky, J.D.M. Wright, M. Youngson, H. Hanche-Olsen, E. St{\o}rmer, E.M. Alfsen and F.W. Shultz, among others.  Detailed definitions and complete sources of references can be found in subsection \ref{subsection background} and in the monographs \cite{AlfseShultz2001, AlfseShultz2003, HOS, Cabrera-Rodriguez-vol1, Cabrera-Rodriguez-vol2}. The notions of invertibility, spectrum and unitaries admit a perfect translation to the setting of unital JB$^*$-algebras, and there is a version of the Russo--Dye theorem established by J.D.M. Wright and M. Youngson (see subsection \ref{subsection background} for more details). In the case of a unital JB$^*$-algebra $M$, the subset, $\mathcal{U} (M),$ of all unitaries in $M$ is not a subgroup of $M$ --actually the Jordan product of two unitary elements need not be, in general, a unitary element even in the case of unital C$^*$-algebras--. So, in the Jordan setting we can no longer speak about the unitary group. \smallskip

In the same way that a celebrated theorem, due to Sakai, identifies von Neumann algebras with those C$^*$-algebras which are dual Banach spaces, in the setting of JB$^*$-algebras, those which are dual Banach spaces are called JBW$^*$-algebras and enjoy additional geometric properties. For example, each JBW$^*$-algebra contains a unit element, and every unitary element in a JBW$^*$-algebra $M$ is of the form $\exp(i h)$, where $h$ lies in the self-adjoint part, $M_{sa},$ of $M$. In \cite{CuPe20b}, the first and fifth authors of this note extended the previously-commented result by Hatori and Moln{\'a}r, and proved that the set of unitaries in a JBW$^*$-algebra is a complete invariant too. More concretely, suppose $\Delta: \mathcal{U} (M)\to \mathcal{U} (N)$ is a surjective isometry, where $M$ and $N$ are two JBW$^*$-algebras. Then there exist a unitary $\omega$ in $N,$ a central projection $p\in N$, and a Jordan $^*$-isomorphism $\Phi:M\to N$ such that $$\begin{aligned}\Delta( u ) &=  U_{\omega^*}\left( p \circ \Phi(u) \right) + U_{\omega^*}\left( (\11b{N}-p) \circ\Phi( u)^*\right)\\
	&=  P_2(U_{\omega^*}( p)) U_{\omega^*}(\Phi(u))  + P_2(U_{\omega^*}(\11b{N} -p)) U_{\omega^*}(\Phi( u^*)),
\end{aligned}$$ for all $u\in \mathcal{U}(M)$. Consequently, $\Delta$ admits a (unique) extension to a surjective real linear isometry from $M$ onto $N$ (see \cite[Theorem 3.9]{CuPe20b}).\smallskip

Since unital C$^*$-algebras are unital JB$^*$-algebras when equipped with the natural Jordan product $a\circ b = \frac12 (a b+ba )$, and the same involution and norm; and when both structures are possible, the notions of unitary coincide, we know that the set of unitaries in a unital JB$^*$-algebra $M$ need not be, in general, connected. The connected component of $\mathcal{U}(M)$ containing the unit element will be also called the \emph{principal component}, and will be denoted by $\mathcal{U}^0 (M)$. \smallskip

Despite nowadays we have excellent and thorough academic monographs covering the theory of JB$^*$-algebras (cf., for example, \cite{AlfseShultz2001, AlfseShultz2003, HOS} and the recent books \cite{Cabrera-Rodriguez-vol1, Cabrera-Rodriguez-vol2}), some geometric aspects have not been explored yet. This is the case of the principal component of the set of unitaries. Section \ref{sec: 2 principal component} is aimed to complete our knowledge on the principal component. The starting point is a result by O. Hatori which asserts that for each unital C$^*$-algebra $A$ we have $$ \begin{aligned}
		\mathcal{U}^0(A) &= \{ e^{i h_1} \cdot \ldots \cdot e^{i h_n} \textbf{1} e^{i h_n} \cdot \ldots \cdot e^{i h_1} : n\in \mathbb{N}, \ h_1, \ldots, h_n\in A_{sa} \}
	\end{aligned}$$ 
is open, closed, and path connected, in the (relative) norm topology on $\mathcal{U}(A)$  (see \cite[Lemma 3.2]{Hatori14}). The advantage of the previous description is that it admits an expression in terms of expressions given by Jordan products. In Theorem \ref{t principal component as product of exp} we prove that for each unital JB$^*$-algebra $M$, the principal component of the set of unitaries in $M$ admits the following description $$\begin{aligned} \mathcal{U}^0(M) &=  M^{-1}_{\textbf{1}}\cap\mathcal{U} (M)\\
			&=\left\lbrace  U_{e^{i h_n}}\cdots U_{e^{i h_1}}(\textbf{1}) \colon n\in \mathbb{N}, \ h_j\in M_{sa} \ \forall\ 1\leq j \leq n  \right\rbrace \\
			&= \left\lbrace  u\in \mathcal{U} (M) : \hbox{ there exists } w\in \mathcal{U}^0(M) \hbox{ with } \|u-w\|<2 \right\rbrace
	\end{aligned},$$ where $M^{-1}_{\textbf{1}}$ stands for the principal component of the set $M^{-1}$ of all invertible elements in $M$. Consequently, $\mathcal{U}^0(M)$ is analytically arcwise connected. A central notion in the theory of Jordan algebras has already appeared in the previous description where the $U_a$ operator associated with an element $a$ in a Jordan algebra $M$ is defined by  $$U_{a} (x) =2 (a\circ x) \circ a  - (a\circ a)\circ x, \ \ (x\in M).$$ When an associative algebra $A$ is regarded as a Jordan algebra with respect to the natural Jordan product we have $U_a (x) = a x a $ for all $a,x\in A$.\smallskip

We should admit that the instability of $\mathcal{U} (M)$ under Jordan products is a handicap to find a description of this set in terms of algebraic properties. A similar instability affects the set of invertible elements in a unital Jordan-Banach algebra. In any case, given two invertible (respectively, unitary) elements $a,b$ in a unital Jordan-Banach algebra (respectively, in a unital JB$^*$-algebra) $M$ the element $U_a(b)$ is invertible  (respectively, unitary). Let $M$ be a unital Jordan-Banach algebra. Following \cite{Pe2020}, we shall say that a subset $\mathfrak{M} \subseteq M^{-1}$ is a \emph{quadratic subset} if $U_{\mathfrak{M}} (\mathfrak{M}) \subseteq \mathfrak{M}$. Let us observe that $\mathcal{U} (M)$ is a self-adjoint subset of $M$ whenever $M$ is a unital JB$^*$-algebra, in such a case we also prove that $\mathcal{U}^{0} (M)$ is the smallest quadratic subset of $\mathcal{U} (M)$ containing the set $e^{i M_{sa}}$ (see Proposition \ref{p self adjoint quadratic subset princ comp}).\smallskip

Surjective isometries between different connected components of the unitary sets of two unital JB$^*$-algebras are studied in section \ref{sec: surj isom}. In a first result we establish that for each surjective isometry $\Delta: \mathcal{U}^0 (M)\to \mathcal{U}^0 (N)$ between the principal components of two unital JB$^*$-algebras, there exist $k_1,\ldots,k_n\in N_{sa}$, a central projection $p\in N$ and a Jordan $^*$-isomorphism $\Phi:M\to N$ such that $$\begin{aligned}\Delta(u) &= p\circ U_{e^{i k_n}} \ldots  U_{e^{i k_1}}  \Phi (u) + (\11b{N}-p)\circ \left( U_{e^{-i k_n}} \ldots  U_{e^{-i k_1}} \Phi(u)\right)^*,
\end{aligned}$$ for all $u\in \mathcal{U}^0(M)$. Consequently, $M$ and $N$ are Jordan $^*$-isomorphic, and there exists a surjective real linear isometry {\rm(}i.e., a real linear triple isomorphism{\rm)} from $M$ onto $N$ whose restriction to $\mathcal{U}^0 (M)$ is $\Delta$ (see Theorem \ref{t surj isom principal components}). Surjective isometries between connected components of the unitary sets of two unital JB$^*$-algebras are described in Theorem \ref{t surj isom secondary components}.\smallskip

Corollary 3.8 in \cite{CuPe20b} proves that two unital JB$^*$-algebras $M$ and $N$ are isometrically isomorphic as (complex) Banach spaces if, and only if, they are isometrically isomorphic as real Banach spaces if, and only if, there exists a surjective isometry between their unitary sets. Corollary \ref{c two unital JB*algebras are isomorphic iff their unitaries are isometric iff their principal components are isometric} in this note improves this conclusion by showing that the following statements are equivalent:
	\begin{enumerate}[$(a)$]
		\item $M$ and $N$ are Jordan $^*$-isomorphic;
		\item There exists a surjective isometry $\Delta: \mathcal{U}(M)\to \mathcal{U}(N)$ satisfying $\Delta (\11b{M})\in \mathcal{U}^0(N)$;
		\item There exists a surjective isometry $\Delta: \mathcal{U}^0(M)\to \mathcal{U}^0(N).$
	\end{enumerate}

The extendibility of a surjective isometry between the unitary sets of two unital JB$^*$-algebras is discussed in section \ref{subsec: extensibility} (see Corollary \ref{c extendibility of surjective isometries} and Proposition \ref{p extra algebraic hypotheses}).\smallskip 

We cannot conclude this section without pointing out the exceptional particularities naturally linked to the category of JB$^*$-algebras. For example, in a unital C$^*$-algebra $A$, any connected component $\mathcal{U}^c (A)$ of $\mathcal{U} (A)$ is isometrically isomorphic to the principal component $\mathcal{U}^0 (A)$ --it suffices to consider the left or right multiplication operator on $A$ by an element in $\mathcal{U}^c (A)$. In the setting of unital JB$^*$-algebras this geometric property is not always true, combining our results with a construction due to R. Braun, W. Kaup and H. Upmeier in \cite{BraKaUp78}, we shall exhibit an example of a unital JB$^*$-algebra $M$ admitting a unitary $w$ such that the connected component of $\mathcal{U} (M)$ containing $w$ is not isometrically isomorphic to the principal component (see Remark \ref{r existence of unitaries with non isomorphic nor isometric connected components in the Jordan setting}). This counterexample reinforces the fact that a unital JB$^*$-algebra admits many different Jordan products and involutions, at least as many as unitaries, some of them produce structures which are not Jordan $^*$-isomorphic. However, by a fascinating result due to W. Kaup, each JB$^*$-algebra admits an essentially unique structure of JB$^*$-triple (see \cite[Proposition 5.5]{Ka83} and subsection \ref{subsection background}). This is essentially the reason for which we employ arguments from JB$^*$-triple theory in some of our arguments.

\subsection{Definitions and background}\label{subsection background} \ \smallskip

Perhaps the idea of having a product which makes stable the self-adjoint part of the space $B(H),$ of all bounded linear operators on a complex Hilbert space $H$, or more generally, the self-adjoint part of a C$^*$-algebra, led mathematicians and physics to consider the natural Jordan product defined by $a\circ b := \frac12 ( a b + ba)$.  The abstract properties of this product are materialized in the notion of Jordan algebra. A complex (respectively,  real) \emph{Jordan algebra} $M$ is a (non-necessarily associative) algebra over the complex (respectively, real) field whose product ``$\circ$'' is commutative and satisfies the so-called \emph{Jordan identity}: \begin{equation}\label{eq Jordan identity} (a \circ b)\circ a^2 = a\circ (b \circ a^2) \hbox{ ($a,b\in M$).} \end{equation}

A \emph{normed Jordan algebra} is a Jordan algebra $M$ equipped with a norm, $\|.\|$, satisfying $\| a\circ b\| \leq \|a\| \ \|b\|$ ($a,b\in M$). A \emph{Jordan-Banach algebra} is a normed Jordan algebra whose norm is complete. Every real or complex associative Banach algebra is a real Jordan-Banach algebra with respect to the natural product. \smallskip

The Jordan identity implies that Jordan algebras are power associative, that is, each subalgebra generated by a single element $a$ in a Jordan algebra $M$ is associative. More concretely, let us set $a^0 = \textbf{1}$ (if $M$ is unital), $a^1 =a$, and $a^{n+1} = a\circ a^n$ ($n\geq 1$). In this case the identity $a^{m+n} = a^{m} \circ a^{n},$ holds for all natural numbers $n,m$ (cf. \cite[Lemma 2.4.5]{HOS}). In the case that $M$ is a Jordan-Banach algebra, the completeness of the norm assures that the series $\displaystyle \exp(a)= e^a =\sum_{n=0}^\infty \frac{a^n}{n!}$ is uniformly convergent on bounded subsets of $M$. Furthermore, the mapping $a\mapsto e^a$ defines an analytic mapping on $M$. The Jordan-Banach subalgebra $C$ of $M$ generated by an element $a$ is a commutative associative subalgebra with respect to the inherited Jordan product \cite[1.1]{Jacobson81}. Suppose $A$ is an associative Banach algebra. For each element $a\in A$, the powers of $A$ with respect to the associative and to the Jordan product coincide, and thus $e^a$ is just the usual exponential in the usual sense for associative Banach algebras.  It follows that $\exp(a)$ has its usual meaning in the Jordan-Banach subalgebra $C$, and, in particular we have \begin{equation}\label{eq exp commute} e^{s a}\circ  e^{t a} = e^{(s + t) a}, \hbox{ for all } s, t \in \mathbb{C}.
\end{equation}

Two elements $a, b$ in a Jordan algebra $M$ are said to \emph{operator commute} if $$(a\circ c)\circ b= a\circ (c\circ b),$$ for all $c\in M$ (cf. \cite[4.2.4]{HOS}). By the \emph{centre} of $M$ 
 we mean the set of all elements of $M$ which operator commute with any other element in $M$. Any element in the centre is called \emph{central}. \smallskip

In a Jordan algebra $M$, the left or right multiplication operator by a fixed element does not play the role played by the left or right multiplication in an associative Banach algebra. The protagonist in the Jordan setting is associated to the $U$ operator. Let $M$ be a Jordan-Banach algebra. Given $a,b\in M$, we shall write $U_{a,b}$ and $M_a$ for the bounded linear operators on $M$ defined by  $$U_{a,b} (x) =(a\circ x) \circ b + (b\circ x)\circ a - (a\circ b)\circ x,\ \hbox{ and } M_{a} (x) =a\circ x$$ for all $x\in M$. The mapping $U_{a,a}$ will be simply denoted by $U_a$. One of the fundamental identities in Jordan theory assures that \begin{equation}\label{eq fundamental identity UaUbUa} U_{U_a(b)}=U_a U_b U_a, \hbox{ for all } a,b \hbox{ in a Jordan algebra } M \end{equation} (see \cite[2.4.18]{HOS}). It can be easily proved, via an induction argument, that the following generalisation of \eqref{eq fundamental identity UaUbUa} holds in every Jordan algebra $M$:
\begin{equation}\label{eq generalised fundamental identity UaUbUa}
U_{U_{a_n}\cdots U_{a_1} (b)}=U_{a_n}\cdots U_{a_1} U_b U_{a_1}\cdots U_{a_n},
\end{equation} for all $n>0$, and all $a_n,\dots, a_1,b\in M$.\smallskip

An element $a$ in a unital Jordan-Banach algebra $M$ is called \emph{invertible} whenever there exists $b\in M$ satisfying $a \circ b = \textbf{1}$ and $a^2 \circ b = a.$ The element $b$ is unique and it will be denoted by $a^{-1}$ (cf. \cite[3.2.9]{HOS} and \cite[Definition 4.1.2]{Cabrera-Rodriguez-vol1}). We know from \cite[Theorem 4.1.3]{Cabrera-Rodriguez-vol1} that an element $a\in M$ is invertible if and only if $U_a$ is a bijective mapping, and in such a case $U_a^{-1} = U_{a^{-1}}$. If an associative Banach algebra $A$ is regarded with its natural Jordan product, the notion of invertibility in the Jordan setting is precisely the usual notion in the associative context. As in the associative setting, the set $M^{-1}$ of all invertible elements in $M$ is open (cf. \cite[Theorem 4.1.7]{Cabrera-Rodriguez-vol1}), and its principal component (i.e. the connected component containing the unit element) will be denoted by $M^{-1}_{\textbf{1}}$.\smallskip

Additional geometric axioms are required to define JB- and JB$^*$-algebras. The definition of JB$^*$-algebras is essentially due to I. Kaplansky, who introduced them as a Jordan generalization of C$^*$-algebras.  A \emph{JB$^*$-algebra} is a complex Jordan-Banach algebra $M$ equipped
with an algebra involution $^*$ satisfying the following geometric axiom: \begin{equation}\label{eq GN axiom JB*} \|U_a (a^*)  \|= \|a\|^3, \ \ (a\in M).
\end{equation} Every C$^*$-algebra $A$ is a JB$^*$-algebra when equipped with its natural Jordan product and the original norm and involution. We observe that in this case, the geometric axiom \eqref{eq GN axiom JB*} writes in the form $\|a a^* a\| = \|a\|^3$ ($a\in A$), which is known to be an equivalent reformulation of the Gelfand-Naimark axiom. M.A. Youngson proved in \cite[Lemma 4]{youngson1978vidav} that the involution of every JB$^*$-algebra is in fact an isometry.\smallskip

A \emph{JB-algebra} is a real Jordan-Banach algebra $J$ in which the norm satisfies
the following two axioms for all $a, b\in J$:\begin{enumerate}[$(i)$]\item $\|a^2\| =\|a\|^2$;
	\item $\|a^2\|\leq \|a^2 +b^2\|.$
\end{enumerate} As pointed out by Kaplansky, the hermitian part, $M_{sa}$, of a JB$^*$-algebra, $M$, is always a JB-algebra. The converse implication was open for some time, and it was shown to be true in a celebrated result due to J.D.M. Wright, asserting that the complexification of every JB-algebra is a JB$^*$-algebra (see \cite{Wri77}). The reader who may feel the necessity of additional details and explanations for any of the basic results on JB-algebras, JB$^*$-algebras and JB$^*$-triples (whose definition appears below) can consult the monographs \cite{HOS,Chu2012, AlfseShultz2001, AlfseShultz2003, Cabrera-Rodriguez-vol1} and \cite{Cabrera-Rodriguez-vol2}. \smallskip

Let $M$ and $N$ be JB$^*$-algebras. A linear mapping $\Phi: M\to N$ is called a \emph{Jordan homomorphism} if $\Phi (a\circ b) = \Phi (a) \circ \Phi (b)$ ($a,b\in M$), and a \emph{Jordan $^*$-homomorphism} if it is a Jordan homomorphism and $\Phi (a)^* = \Phi(a^*)$ for all $a\in M$. A Jordan $^*$-isomorphism is a Jordan $^*$-homomorphism which is also a bijection. A real linear mapping from $M$ to $N$ preserving Jordan products will be called a real linear Jordan homomorphism. Real linear Jordan $^*$-homomorphisms are similarly defined. A mapping between unital JB$^*$-algebras is called \emph{unital} if it sends the unit in the domain to the unit in the codomain.\smallskip

A JBW$^*$-algebra is a JB$^*$-algebra which is also a dual Banach space. Norm-closed Jordan $^*$-subalgebras of C$^*$-algebras are called \emph{JC$^*$-algebras}. JC$^*$-algebras which are also dual Banach spaces are called \emph{JW$^*$-algebras}. Any JW$^*$-algebra is a weak$^*$-closed Jordan $^*$-subalgebra of a von Neumann algebra. The reader should be warned about the existence of exceptional JB$^*$-algebras which cannot be embedded as Jordan $^*$-subalgebras of some $B(H)$ (see \cite[Corollary 2.8.5]{HOS}, \cite[Example 3.1.56]{Cabrera-Rodriguez-vol1}).\smallskip

We recall that an element $u$ in a unital JB$^*$-algebra $M$ is a \emph{unitary} if it is invertible and its inverse coincides with $u^*$. As in the associative setting, we shall denote by $\mathcal{U}(M)$ the set of all unitary elements in $M$. It is known that the identities $$ U_u(u^*)=u, \hbox{ and } U_u((u^*)^2)=\textbf{1},$$ hold for every unitary element in $M$ {\rm(}cf. \cite[Theorem 4.1.3]{Cabrera-Rodriguez-vol1}{\rm)}. Let us observe that if a unital C$^*$-algebra is regarded as a JB$^*$-algebra both notions of unitaries coincide. An element $s$ in a unital JB-algebra $J$ is called a \emph{symmetry} if $s^2 =\11b{J}$. If $M$ is a JB$^*$-algebra, the symmetries in $M$ are defined as the symmetries in its self-adjoint part $M_{sa}$.\smallskip

In a JBW$^*$-algebra $M$ the set $\mathcal{U}(M)$ is path connected and coincides with the set $\{e^{i h} : h\in M_{sa}\}$ (cf. \cite[Remark 3.2$(9)$]{CuPe20b}). However, as in the case of unital C$^*$-algebras, in a general unital JB$^*$-algebra the set of unitaries might contain many different connected components. The principal component of $\mathcal{U}(M)$ is the connected component containing the unit, and it will be denoted by $\mathcal{U}^0(M)$.\smallskip

\begin{remark}\label{r isotopes remark} A unital JB$^*$-algebra admits, at least, as many different products as unitaries it contains. More concretely, let $M$ be a unital JB$^*$-algebra with unit $\textbf{1}$. To illustrate the statement, let us recall several basic properties.
\begin{enumerate}[$(a)$]\item For each unitary $u\in M$, the Banach space underlying $M$ becomes a unital JB$^*$-algebra with unit $u$ for the Jordan product and involution defined by $$ x\circ_u y :=U_{x,y}(u^*) \hbox{ and } x^{*_u} :=U_u(x^*), \hbox{ respectively.}$$ This new unital JB$^*$-algebra $(M,\circ_u,*_u)$ is called the \emph{$u$-isotope} of $M$, and will be denoted by $M(u)$ {\rm(}see \cite[Lemma 4.2.41$(i)$]{Cabrera-Rodriguez-vol1}{\rm)}. For each $a\in M$, we shall write $a^{(n,u)}$ for the $n$-th power of $a$ in the JB$^*$-algebra $M(u)$;
\item Fix a unitary $u\in M$. The JB$^*$-algebras $M$ and $(M,\circ_u,*_u)$ have the same sets of invertible and unitary elements, that is, $M^{-1} = (M(u))^{-1}$ and  $\mathcal{U} (M) =  \mathcal{U} (M(u))$  {\rm(}see \cite[Lemma 4.2$(ii)$ and Theorem 4.6]{Sidd2007TopologicalRank1} or \cite[Lemma 4.2.41$(ii)$]{Cabrera-Rodriguez-vol1}{\rm)};
\item For each unitary $u\in M$, the mapping $U_u : M\to M$ is a surjective linear isometry. Moreover, when $U_u$ is regarded as an operator between the JB$^*$-algebras $ (M,\circ_{u^*},*_{u^*})$ and $(M,\circ_u,*_u)$, it is a Jordan $^*$-isomorphism {\rm(}see \cite[Theorem 4.2.28]{Cabrera-Rodriguez-vol1}, \cite{WriYou78}{\rm)}.
\end{enumerate}
\end{remark}

Beside the deep connections with holomorphic theory on arbitrary complex Banach spaces and the classification of bounded symmetric domains in this complex setting which led W. Kaup to introduce those complex Banach spaces called JB$^*$-triples in \cite{Ka83}, for the purposes of this note, we shall employ the JB$^*$-triple structure associated with each JB$^*$-algebra to simplify all possible Jordan products and involutions under a triple product which encompasses all of them. A \emph{JB$^*$-triple} is a complex Banach space $E$ equipped with a continuous triple product $\J ... : E\times E\times E \to E,$ $(a,b,c)\mapsto \{a,b,c\},$ which is bilinear and symmetric in $(a,c)$ and conjugate linear in $b$, and satisfies the following axioms for all $a,b,x,y\in E$:
\begin{enumerate}[{\rm (a)}] \item $L(a,b) L(x,y) = L(x,y) L(a,b) + L(L(a,b)x,y)
	- L(x,L(b,a)y),$ where $L(a,b):E\to E$ is the operator defined by $L(a,b) x = \J abx;$
	\item $L(a,a)$ is a hermitian operator with non-negative spectrum;
	\item $\|\{a,a,a\}\| = \|a\|^3$.\end{enumerate}\smallskip

Examples of JB$^*$-triples include all C$^*$-algebras and JB$^*$-algebras with the triple products of the form \begin{equation}\label{eq product operators} \J xyz =\frac12 (x y^* z +z y^* x),\end{equation}  and \begin{equation}\label{eq product jordan}\J xyz = (x\circ y^*) \circ z + (z\circ y^*)\circ x -
(x\circ z)\circ y^*, \end{equation} respectively (\cite[Fact 4.1.41 and Theorem 4.1.45]{Cabrera-Rodriguez-vol1}). Actually, when a C$^*$-algebra is regarded as a JC$^*$-algebra with its natural Jordan product, the expression in \eqref{eq product jordan} coincides with the triple product in \eqref{eq product operators}. The product in \eqref{eq product operators} is actually valid to define a structure of JB$^*$-triple on spaces $B(H,K)$ of bounded linear operators between complex Hilbert spaces (infinite dimensional examples of rectangular matrices), complex Hilbert spaces, and J$^*$-algebras in the sense introduced by  L. Harris in \cite{Harris74,Harris81}--i.e. closed complex-linear subspaces of $B(H,K)$ which are closed for the triple product in \eqref{eq product operators}--.\smallskip

The triple product of every JB$^*$-triple is a non-expansive mapping, that is,
\begin{equation}\label{eq triple product non-expansive} \|\{a,b,c\}\|\leq \|a\| \|b\| \|c\|\ \hbox{ for all } a,b,c \hbox{ (see \cite[Corollary 3]{FriRu86}).}
\end{equation}

A (real linear) triple homomorphism between JB$^*$-triples $E$ and $F$ is a (real) linear mapping $\Phi : E\to F$ preserving triple products, i.e., $$T\{a,b,c\} = \{T(a),T(b),T(c)\}, \hbox{ for all $a,b,c\in E$.}$$ It is well known that every unital triple homomorphism between unital JB$^*$-algebras is a Jordan $^*$-homomorphism.\label{eq unital triple hom are Jordan star hom}\smallskip

Among the amazing geometric properties of those complex Banach spaces in the class of JB$^*$-triples we find a Banach-Stone type theorem, proved by W. Kaup, asserting that a linear bijection $T$ between two JB$^*$-triples is an isometry if and only if it is a triple isomorphism (cf. \cite[Proposition 5.5]{Ka83}, see also \cite[Theorem 2.2.28]{Cabrera-Rodriguez-vol1} or \cite{DaFriRu} and \cite{FerMarPe2004} for alternative proofs). Thanks to this result we can understand now that, despite of the existence of many different Jordan products and involutions on a unital JB$^*$-algebra, they all produce the same triple structure. Namely, if $u_1$ and $u_2$ are two unitaries in a unital JB$^*$-algebra $M$, we know that $M$ is a JB$^*$-triple for the following three triple products: $$\begin{aligned}\{x,y,z\} &= (x\circ y^*) \circ z + (z\circ y^*)\circ x -
		(x\circ z)\circ y^* = U_{x,z} (y^*), \\
\{x,y,z\}_1 &= (x\circ_{u_1} y^{*_{u_1}}) \circ_{u_1} z + (z\circ_{u_1} y^{*_{u_1}})\circ_{u_1} x -
		(x\circ_{u_1} z)\circ_{u_1} y^{*_{u_1}}, \hbox{ and } \\
\{x,y,z\}_2 &= (x\circ_{u_2} y^{*_{u_2}}) \circ_{u_2} z + (z\circ_{u_2} y^{*_{u_2}})\circ_{u_2} x -
		(x\circ_{u_2} z)\circ_{u_2} y^{*_{u_2}}.
\end{aligned}$$ Since the identity mapping is a surjective linear isometry between any of these three structures, the three triple products given above must coincide.\label{eq coincidence of the Jordan triple products}\smallskip

It should be remarked that Kaup's Banach--Stone theorem is not, in general, true for surjective real linear isometries between JB$^*$-triples (cf. \cite[Remark 2.7]{Da}, \cite[\S 4]{IsKaRo95} and \cite{FerMarPe}). However, surjective real linear isometries between C$^*$-algebras and JB$^*$-algebras are all triple isomorphisms (see \cite[Corollaries 3.2 and 3.3]{Da} and \cite[Corollary 3.4]{FerMarPe}). We further know that every unital surjective real linear isometry between unital JB$^*$-algebras is a real linear Jordan $^*$-isomorphism \cite[Corollary 3.2]{Da}.\smallskip

For each couple of elements $a,b$ in a JB$^*$-triple $E$, we denote by $Q(a,b)$ the conjugate linear operator on $E$ defined by $Q(a,b)(x) = \{a,x,b\}$. We shall write $Q(a)$ for $Q(a,a)$. In the case of a JB$^*$-algebra $M$, the $Q$ operator is intrinsically related to the $U$ operator by the identity $$Q(a,b) (x) = \{a,x,b\} = U_{a,b} (x^*), \ \hbox{ for all } x\in M.$$

We have already mentioned that the involution of every unital JB$^*$-algebra $M$ is an isometry. We can actually apply the previous Kaup's Banach-Stone theorem to the involution as a surjective linear isometry from $M$ onto the Banach space $M$ equipped with the same norm but replacing the product by scalars by the one given by $\lambda \odot x = \overline{\lambda} x$, or it can be directly checked that the involution is a (conjugate linear) triple isomorphism. Consequently, for any $a$ and $b$ in $M$, we have that \begin{equation}\label{eq simple U-op adjoint}(U_a(b))^*=(\J{a}{b^*}{a})^*=\J{a^*}{b}{a^*}=U_{a^*}(b^*).\end{equation} By an induction argument we have
\begin{equation}\label{eq U-op adjoint}
(U_{a_n}\cdots U_{a_1}(a_0))^*=U_{a_n^*}\cdots U_{a_1^*}(a_0^*),\end{equation}
for all natural $n$ and $a_n,\dots, a_1, a_0\in M$.\smallskip

Each element $e$ in a JB$^*$-triple $E$ satisfying $\J eee=e$ is called a \emph{tripotent}. If a C$^*$-algebra $A$ is regarded as a JB$^*$-triple with the triple product in \eqref{eq product operators}, the partial isometries in $A$ are precisely the tripotents of $A$. Each tripotent $e\in E$ determines a \emph{Peirce decomposition}\label{eq Peirce decomposition} of $E$ in the form
$$E= E_{2} (e) \oplus E_{1} (e) \oplus E_0 (e),$$ where $E_j (e)=\{ x\in E : \J eex = \frac{j}{2}x \}$ ($j=0,1,2$) is called the \emph{Peirce $j$-subspace}.\smallskip

Triple products among elements in Peirce subspaces satisfy the following \emph{Peirce arithmetic}: $$\begin{aligned} \J {E_{i}(e)}{E_{j} (e)}{E_{k} (e)} &\subseteq E_{i-j+k} (e)\  \hbox{ if $i-j+k \in \{ 0,1,2\},$ }\\
\J {E_{i}(e)}{E_{j} (e)}{E_{k} (e)} & =\{0\} \hbox{ if $i-j+k \notin \{ 0,1,2\},$ }
\end{aligned}$$ and $\J {E_{2} (e)}{E_{0}(e)}{E} = \J {E_{0} (e)}{E_{2}(e)}{E} =\{0\}.$ Consequently, each Peirce subspace $E_j(e)$ is a JB$^*$-subtriple of $E$.\smallskip

The natural projection, $P_{k_{}}(e),$ of $E$ onto $E_{k} (e)$ is called the \emph{Peirce $k$-projection}. It is known that $P_{2}(e) = Q(e)^2,$ $P_{1}(e) =2(L(e,e)-Q(e)^2),$ and $P_{0}(e) =\hbox{Id}_E - 2 L(e,e) + Q(e)^2.$ Furthermore, $\|P_{k_{}}(e)\|\leq 1$ for all $k=0,1,2$ (cf. \cite[Corollary 1.2]{FriRu85}).\smallskip

It is worth remarking that if $e$ is a tripotent in a JB$^*$-triple $E$, the Peirce 2-subspace $E_2 (e)$ is a unital JB$^*$-algebra with unit $e$,
product $x\circ_e y := \J xey$ and involution $x^{*_e} := \J exe$, respectively (cf. \cite[Theorem 4.1.55]{Cabrera-Rodriguez-vol1}).\label{eq Peirce-2 is a JB-algebra} 
\smallskip

The adjective ``unitary'' is also employed to denote some special tripotents. More precisely, a tripotent $e$ in a JB$^*$-triple $E$ is called \emph{unitary} 
if $E_2 (e) = E$. 
However, there is no risk of controversy, since, as shown by R. Braun, W. Kaup and H. Upmeier in \cite[Proposition 4.3]{BraKaUp78}, the unitaries in a unital JB$^*$-algebra $M$ are precisely the unitary tripotents in $M$ when the latter is regarded as a JB$^*$-triple. We can therefore conclude that every surjective real linear isometry $T$ between two unital JB$^*$-algebras $M$ and $N$ maps $\mathcal{U} (M)$ onto $\mathcal{U} (N)$.  In particular, for each unitary $u\in M$, the mapping $U_u : M\to M$ is a triple isomorphism satisfying $U_u \left( \mathcal{U}(M) \right) = \mathcal{U}(M)$.\smallskip


For later purposes we conclude this subsection on background results and definitions by recalling a Jordan version of the Stone's one-parameter theorem for norm continuous one parameter quadratic maps which has been borrowed from the recent paper \cite{CuPe20b}.

\begin{theorem}\label{t Jordan unitary groups version of Stone's theorem}\cite[Theorem 3.1]{CuPe20b} Let $M$ be a unital JB$^*$-algebra. Suppose $\{u(t):t\in \mathbb{R}\}$ is a family in $\mathcal{U} (M)$ satisfying $u(0) =\textbf{1},$ and $U_{u(t)} (u(s)) = u(2 t +s),$ for all $t,s\in \mathbb{R}$. We also assume that the mapping $t\mapsto u(t)$ is norm continuous. Then there exists $h\in M_{sa}$ such that $u(t) = e^{i t h}$ for all $t\in \mathbb{R}$.
\end{theorem}

\section{Range tripotents and the principal component of the unitary set in a JB$^*$-algebra}\label{sec: 2 principal component}

The principal components of the sets of invertible elements, $A^{-1},$ and of the unitary group, $\mathcal{U}(A),$ of a unital C$^*$-algebra $A$ are the connected components of these subsets containing the unit $\textbf{1}\in A$. These principal components will be denoted by $A_{\textbf{1}}^{-1}$ and $\mathcal{U}^0(A)$, respectively. These two sets are well known, studied and described in the literature. More concretely, for any (associative) unital real or complex Banach algebra $A$, $A^{-1}$ is an open subset of $A$ which is stable under associative products. The principal component of $A^{-1}$ is precisely the least subgroup of $A^{-1}$ containing $\exp(A)=\{ e^{a} : a\in A\},$ it is path connected, and it can be algebraically described in the following terms: $$ A_{\textbf{1}}^{-1} = \{ e^{a_1} \cdot \ldots \cdot e^{a_n} : n\in \mathbb{N}, \ a_1, \ldots, a_n\in A\}$$
(see \cite[Propositions 8.6 and 8.7]{BonsDun73}, \cite[page 12]{Tak} and \cite[Theorem 2.3.1]{Li2003}).\smallskip

The principal component of the unitary group of a unital C$^*$-algebra $A$ has been also studied and described by different authors. It is known that
\begin{equation}\label{eq algebraic descirption of the principal component for unital Cstar algebras} \begin{aligned}
\mathcal{U}^0(A) &= \{ e^{i h_1} \cdot \ldots \cdot e^{i h_n} : n\in \mathbb{N}, \ h_1, \ldots, h_n\in A_{sa} \}\\ &= \{ e^{i h_1} \cdot \ldots \cdot e^{i h_n} \textbf{1} e^{i h_n} \cdot \ldots \cdot e^{i h_1} : n\in \mathbb{N}, \ h_1, \ldots, h_n\in A_{sa} \}
\end{aligned}
 \end{equation} is open, closed, and path connected, in the (relative) norm topology on $\mathcal{U}(A)$ --the first equality in \eqref{eq algebraic descirption of the principal component for unital Cstar algebras} can be found in Kadison--Ringrose book \cite[Exercises 4.6.6 and 4.6.7]{KR1}, while the second equality, which is closer to the Jordan structure, is due to O. Hatori \cite[Lemma 3.2]{Hatori14}. Another interesting result asserts that, by the continuity of the module mapping $x\mapsto |x|:=(x^* x)^\frac12$, we have $$A^{-1}_{\textbf{1}} \cap \mathcal{U}(A) =\mathcal{U}^0(A) \ \ \ \hbox{{\rm(}see \cite[Exercise 2, page 56]{Tak}{\rm).}}$$ The reader should be warned that ``polar decompositions'' are not available for elements in a JB$^*$-algebra, so this kind of arguments are not valid in the Jordan setting.\smallskip

To the best of our knowledge, in the setting of real C$^*$-algebras a full description of the principal component of the unitary group remains unknown. For example $\mathcal{U} (\mathbb{R}) = \{-1,1\}$ has two connected components, $\mathcal{U}^0 (\mathbb{R}) = \{1\} = \{e^0\}$ and $\mathbb{R}_{skew} = \{0\}$, while $\mathcal{U} (\mathbb{C}) = \mathbb{T} =\mathcal{U}^0 (\mathbb{C}) = \{e^k : k \in\mathbb{C}_{skew}\}$  is connected.\smallskip

In the case of a unital C$^*$-algebra $A$, the sets $A^{-1}$, $A_{\textbf{1}}^{-1},$ $\mathcal{U}(A)$, $\mathcal{U}^0(A)$ are closed for the associative product but they are not, in general, stable under Jordan products.\smallskip

In the case of Jordan-Banach algebras, B. Aupetit posed the question whether the principal component of the set of invertible elements in a Jordan-Banach algebra is analytically arcwise connected in \cite[comments after the proof of Theorem 3.1]{Aupetit95}. A complete positive solution to this question was found by O. Loos in \cite{Loos96} (see also \cite[Theorem 4.1.111]{Cabrera-Rodriguez-vol1}), where the following conclusion was established:

\begin{theorem}\label{t Loos principal component invertible}{\rm\cite[Corollary]{Loos96}} Let $M$ be a unital Jordan-Banach algebra. Then the principal component of $M^{-1}$ is the set $$M_{\textbf{1}}^{-1} = \left\{ U_{\exp(a_1)}\cdots U_{\exp(a_n)} (\textbf{1}) : n\in \mathbb{N},\ a_1, \ldots, a_n\in M \right\},$$ in particular, $M_{\textbf{1}}^{-1}$ is open and closed in $M^{-1},$ and each connected component of $M^{-1}$ is analytically arcwise connected.
\end{theorem}

However, as long as we know the description of the principal component of the set of unitaries in a unital JB$^*$-algebras has not been explicitly treated in the literature. By employing some tools in the recent note \cite{CuPe20b} we shall complete the state-of-the-art in the Jordan setting. One of the main geometric properties derived from the local Gelfand theory in C$^*$-algebras assures that for each unitary $u$ in a unital C$^*$-algebra $A$ with $\|\textbf{1} -u\|<2$ we can always find a hermitian element $h\in A_{sa}$ satisfying $u = e^{ih}$ (see \cite[Exercise 4.6.6]{KR1}). Building upon the celebrated \emph{Shirshov-Cohn theorem}, which proves that the JB$^*$-subalgebra of a JB$^*$-algebra generated by two self-adjoint elements (and the unit element) is a JC$^*$-algebra, that is, a JB$^*$-subalgebra of some $B(H)$ (cf. \cite[Theorem 7.2.5]{HOS} and \cite[Corollary 2.2]{Wri77}), this fact was extended to the setting of unital JB$^*$-algebras in \cite{CuPe20b}.

\begin{lemma}\label{l untaries at distance <2 in a unital JB$^*$-algebra}\cite[Lemma 2.2]{CuPe20b} Let $u,v$ be two unitaries in a unital JB$^*$-algebra $M$. Let us suppose that $\|u-v\| <2.$ Then there exists a self-adjoint element $h$ in the $u$-isotope JB$^*$-algebra $M(u)= (M,\circ_u,*_u)$ such that $$v = e_u^{ih}=\exp_u (i h) = \sum_{n=0}^{\infty} \frac{i^n}{n!} h^{(n,u)},$$ that is, the exponential is computed in the JB$^*$-algebra $M(u)$. Consequently, if $\|\textbf{1}-v\| <2,$ there exists a hermitian element $h\in M_{sa}$ satisfying $v = e^{i h}$ and $U_{e^{i \frac{h}{2}}} (\textbf{1}) = v$.
\end{lemma}

Let $M$ be a unital JB$^*$-algebra. As in the setting of C$^*$-algebras, the symbol $\mathcal{U}^0(M)$ will stand for the principal component of the set of all unitaries in $M$, that is, the connected component of $\mathcal{U}(M)$ containing the unit $\11b{M}$ of $M$.\smallskip

Let us consider a unitary element $u$ in a unital JB$^*$-algebra $M$, and the $u$-isotope $M(u):=(M,\circ_u,*_u)$ as defined in subsection \ref{subsection background}. It is known that $\mathcal{U}(M) =\mathcal{U}(M(u))$ (cf. Remark \ref{r isotopes remark}$(b)$). Let us observe that the  unital JB$^*$-algebras $M$ and $M(u)$ share the same underlying Banach spaces, the quoted unitary sets are equipped with the same norm and distance. Let us suppose that $u$ lies in $\mathcal{U}^0(M)$. Since the connected components of a topological space are disjoint, we can affirm that
\begin{equation}\label{eq princ components of isotopes coincide}
\mathcal{U}^0(M(u))=\mathcal{U}^0(M), \text{ for every $u\in \mathcal{U}^0(M)$.}
\end{equation}

Section \ref{sec: surj isom} below will be devoted to study and describe the surjective isometries between the principal components of the unitary sets of two unital JB$^*$-algebras. For this purpose, we shall first establish an algebraic characterization of the set of unitary elements in terms of $U$ operators of exponentials of skew symmetric elements, like the one commented in the case of unital C$^*$-algebras (cf. \eqref{eq algebraic descirption of the principal component for unital Cstar algebras}).\smallskip

Let us recall the local Gelfand theory to the readers. In a C$^*$-algebra, the C$^*$-subalgebra $A$ generated by a non-normal element cannot be identified, via Gelfand theory, with a commutative C$^*$-algebra of the form $C_0(L)$, for a locally compact Hausdorff space $L$. However, if we regard $A$ as a JB$^*$-triple, the JB$^*$-subtriple generated by a single element is always representable as a $C_0(L)$ space. Let $0\neq a$ be an element in a JB$^*$-triple $E$. The JB$^*$-subtriple of $E$ generated by $a$ (and denoted by $E_a$) coincides with the norm closure of the linear span of the odd powers of $a$ defined as $a^{[1]}= a$, $a^{[3]} = \J aaa$, and $a^{[2n+1]} := \J aa{a^{[2n-1]}},$ $(n\in \NN)$. The local Gelfand theory for JB$^*$-triples assures that $E_a$ is isometrically JB$^*$-triple isomorphic to some $C_0 (\Omega_{a})$ for a unique compact Hausdorff space $\Omega_{a}$ contained in the set $[0,\|a\|],$ such that $0$ cannot be an isolated point in $\Omega_a$, where, along this note, the symbol $C_0 (\Omega_{a})$ will stand for the Banach space of all complex-valued continuous functions on $\Omega_a$ vanishing at $0$ if $0\in \Omega_a$, and $C (\Omega_{a})$ otherwise. It is further known that there exists a triple isomorphism $\Psi : E_a \to C_{0}(\Omega_a)$ satisfying $\Psi (a) (t) = t$ $(t\in \Omega_a)$ (cf. \cite[Corollary 4.8]{Ka0}, \cite[Corollary 1.15]{Ka83} and \cite{FriRu85}, see also \cite[Lemma 3.2, Corollary 3.4 and Proposition 3.5]{Ka96}, \cite[Theorem 4.2.9]{Cabrera-Rodriguez-vol1} or \cite[Theorem 3.1.12]{Chu2012}). The set $\Omega_a$ is called \emph{the triple spectrum} of $a$ (in $E$), and it does not change when computed with respect to any JB$^*$-subtriple $F\subseteq E$ containing $a$ \cite[Proposition 3.5]{Ka96}. 
\smallskip

As in the associative setting, the local Gelfand theory opens the door to apply a \emph{continuous triple functional calculus}. Given an element $a$ in a JB$^*$-triple $E$, let $\Omega_a$ and $\Psi : E_a \to C_{0}(\Omega_a)$ denote the triple spectrum of $a$ and the triple isomorphism given in previous paragraph. For each continuous function $f\in C_{0}(\Omega_a)$, we set $f_t(a) :=\Psi^{-1} (f)$, and we call it the \emph{continuous triple functional calculus} of $f$ at $a$. If $p(\lambda) = \alpha_1 \lambda + \alpha_3 \lambda^3 + \ldots + \alpha_{2 n -1} \lambda^{2 n -1}$ is an odd polynomial with complex coefficients, it is easy to check that $p_t (a) = \alpha_1 a + \alpha_3 a^{[3]} + \ldots + \alpha_{2 n -1} a^{[2 n -1]}$.\smallskip

Another example, the function $g(t) =\sqrt[3]{t}$ produces $g_t(a) = a^{[\frac{1}{3}]}$. Let us note that $g_t(a) = a^{[\frac{1}{3}]}$ is the unique \emph{cubic root} of $a$ in $E_a$, i.e. the unique element $a^{[\frac13 ]}\in E_a$ satisfying \begin{equation}\label{eq existence of cubic roots} \J {a^{[\frac13 ]}}{a^{[\frac13 ]}}{a^{[\frac13 ]}}=a.
\end{equation} The sequence $(a^{[\frac{1}{3^n}]})_n$ can be recursively defined by $a^{[\frac{1}{3^{n+1}}]} = \left(a^{[\frac{1}{3^{n}}]}\right)^{[\frac 13]}$, $n\in \NN$. It is easy to see that the sequence $(a^{[\frac{1}{3^n}]})_n$ need not be, in general, norm convergent in $E$ (consider, for example, the element $a(t) =t$ in $E=C([0,1])$).\smallskip

The functional analysis provides another topology to assure the convergence of the above sequence. Let us recall that a \emph{JBW$^*$-triple} is a JB$^*$-triple which is also a dual Banach space (with a unique isometric predual \cite{BarTi86}). The triple product of every JBW$^*$-triple is separately weak$^*$-continuous (cf. \cite{BarTi86}), and the second dual, $E^{**}$, of a JB$^*$-triple $E$ is always a JBW$^*$-triple under a triple product which extends the original triple product in $E$ (cf. \cite{Di86}).\smallskip

If $W$ is a JBW$^*$-triple (in particular when we regard $E$ as a JB$^*$-subtriple of $E^{**}$), the sequence $(a^{[\frac{1}{3^n}]})_n$ converges in the weak$^*$-topology of $W$ to a (unique) tripotent denoted by $r(a)$. The tripotent $r(a)$ is called the \emph{range tripotent} of $a$ in $W$. This range tripotent $r(a)$ can be characterized as the smallest tripotent $e\in W$ satisfying that $a$ is positive in the JBW$^*$-algebra $W_{2} (e)$ (compare \cite[Lemma 3.3]{EdRu96}).\smallskip

Given an element $a$ in a JB$^*$-triple $E$, we shall denote by $r(a)$ the range tripotent of $a$ in $E^{**}$. There exist examples in which we can guarantee that the element $r(a)$ lies in $E$. According to the usual notation followed in \cite{BurKaMoPeRa,FerGarSanSi,Ka96,JamPerSidTah2015}, an element $a$ in a JB$^*$-triple $E$ is called \emph{von Neumann regular} if $a\in Q(a) (E)$; if $a\in Q(a)^2 (E)$ we say that $a$ is \emph{strongly von Neumann regular}. For each von Neumann regular element $a$ in a JB$^*$-triple $E$ there might exist many different elements $c$ in $E$ satisfying $Q(a)(c) =a$. However, it is known (cf. \cite[Theorem 1]{FerGarSanSi} and
\cite[Lemma 4.1, Lemma 3.2 and comments after its proof]{Ka96}, \cite[Theorem 2.3 and Corollary 2.4]{BurKaMoPeRa}) that the following statements are equivalent:
\begin{enumerate}[$(a)$]\item $a$ is von Neumann regular;
\item $a$ is strongly von Neumann regular;
\item There exists $b\in E$ such that $Q(a) (b) =a,$ $Q(b) (a) =b$ and $[Q(a),Q(b)]:=Q(a)\,Q(b) - Q(b)\, Q(a)=0$;
\item $0\notin \Omega_a,$ \emph{the triple spectrum} of $a$;
\item $Q(a)$ has norm-closed range;
\item The range tripotent $r(a)$ of $a$ lies in $E$ and $a$ is positive and invertible in the JB$^*$-algebra $E_2 (r(a))$.
\end{enumerate} The element $b$ appearing in statement $(c)$ above is unique, and it will be called the \emph{generalized inverse} of $a$ in $E$ (denoted by $a^{\dag}$). Henceforth we shall write Reg$_{_{vN}}(E)$ for the set of all von Neumann regular elements in $E$. It is further known that $L(a,a^{\dag}) = L(a^{\dag},a) = L(r(a),r(a))$ and $Q(a) Q(a^{\dag}) = Q(a^{\dag}) Q(a) = P_2(r(a))$ (see, for example, \cite[page 589]{JamPerSidTah2015}, \cite[comments in page 192]{BurKaMoPeRa} or \cite[Lemma 3.2 and subsequent comments]{Ka2001}).\smallskip

If $a$ is an invertible element in a unital JB$^*$-algebra $M$ and the latter is regarded as a JB$^*$-triple, then $a$ is von Neumann regular and its range tripotent, $r(a),$ is a unitary element in $M$, that is, \begin{equation}\label{eq range of invertible is unitary} a\in M^{-1} \Rightarrow r(a)\in \mathcal{U} (M) \ \ \hbox{ (cf. \cite[Lemma 2.2 and Remark 2.3]{JamPerSidTah2015}).}
\end{equation} Namely, in this case, by denoting $z^{-1}$ the inverse of $z$ in $M$ we have $$\begin{aligned}Q(z) ((z^{-1})^*) &= \{z, (z^{-1})^*, z\} = U_{z} (z^{-1}) = z; \\
Q((z^{-1})^*) (z) &= \{(z^{-1})^*,z,(z^{-1})^*\} = U_{(z^{-1})^*} (z^*) = (z^{-1})^*;\\
 Q(z) Q((z^{-1})^*) (x) &=  U_z \left( U_{(z^{-1})^*} (x^*) \right)^* = U_z U_{z^{-1}} (x) =  U_z U_{z}^{-1} (x) = x; \\
 Q((z^{-1})^*) Q(z) (x) &= U_{(z^{-1})^*} \left( U_z (x^*) \right)^* = U_{(z^*)^{-1}} U_{z^*} (x) = U_{z^*}^{-1} U_{z^*} (x) = x,
\end{aligned}$$ for all $x\in M$, which proves that $z$ is von Neumann regular with $z^{\dag} = (z^{-1})^* = (z^*)^{-1},$ and since $P_2(r(z)) = Q(z) Q(z^{\dag}) = Id$, the tripotent $r(z)$ is a unitary in $M$.\smallskip

We explore next certain properties of the continuous triple functional calculus which might result surprising.  Let us fix an element $a$ in a JB$^*$-triple $E$. Let us consider the isometric triple isomorphism $\Psi: E_a\to C_0(\Omega_a)$ satisfying $\Psi (a) (t) =t $ ($t\in \Omega_a$). The function $g(t) = t^2$ lies in $C_0(\Omega_a)$, and hence $g_t(a) = a^{[2]}$ is an element in $E_a$. The reader should be warned that $a^{[2]}$ is only a square for a local binary product depending on the element $a$. Suppose, additionally, that $a$ is von Neumann regular, that is, $\Omega_a\subset \mathbb{R}^+$. In this case the function $h(t) = \frac1t$ lies in $C_0(\Omega_a)$ and $a^{[-1]} = h_t(a)\in E_a$. It is easy to see that $$\{a,a, a^{[-1]}\} = a, \{a,r(a), a^{[-1]}\} = r(a),\hbox{ and } \{r(a), a^{[-1]}, r(a)\} = a^{[-1]}.$$
The elements $a, a^{[-1]},$ $a^{[2]}$ and $a^{\dag}$ belong to the unital JB$^*$-algebra $E_2 (r(a))$. 
 Having in mind that, by Kaup's theorem, when the triple product of $E$  is restricted to $E_2(r(a))$ it coincides with
$$\{x,y,z\} = (x\circ_{r(a)} y^{*_{r(a)}}) \circ_{r(a)} z + (z\circ_{r(a)} y^{*_{r(a)}})\circ_{r(a)} x -
		(x\circ_{r(a)} z)\circ_{r(a)} y^{*_{r(a)}}, $$ for all $x,y,z\in E_2(r(a)).$ Then $$r(a) = \{a^{\dag}, a, r(a)\} = L(a^{\dag},a) (r(a)) = a^{\dag}\circ_{r(a)} a,$$ $$\begin{aligned}\hbox{ and } a &= L(a^{\dag},a) (a)  = \{a^{\dag}, a, a\} \\
&= (a^{\dag}\circ_{r(a)} a^{*_{r(a)}}) \circ_{r(a)} a + (a\circ_{r(a)} a^{*_{r(a)}})\circ_{r(a)} a^{\dag} - (a^{\dag}\circ_{r(a)} a)\circ_{r(a)} a^{*_{r(a)}}\\
&= (a\circ_{r(a)} a)\circ_{r(a)}  a^{\dag}.\end{aligned}$$ Therefore $ a^{\dag} = a^{[-1]}$.\smallskip

We stop next in a version of \cite[Proposition I.4.10]{Tak}. The proof can provide to the readers a better understanding of the continuous triple functional calculus.

\begin{proposition}\label{p continuity of the square} Let $\Omega\subset \mathbb{R}_0^+$ be a compact set, and let $E$ be a JB$^*$-triple. Suppose  $\mathcal{E}_{\Omega}$ denotes the set of all elements $a\in E$ with $\Omega_a\subseteq \Omega$. If $f:\Omega\to \mathbb{C}$ is a continuous function, then the continuous triple functional calculus: $a\in \mathcal{E}_{\Omega}\mapsto f_t(a)\in E$ is {\rm(}norm{\rm)} continuous. Consequently, the mapping $a\in E\mapsto a^{[2]}\in E$ is continuous.
\end{proposition}

\begin{proof} A variant of the Stone-Weierstrass theorem established in \cite{Ka83} and \cite[Lemma 4.2.8]{Cabrera-Rodriguez-vol1} assures that if $\iota: \Omega\to\mathbb{C}$ denotes the inclusion mapping, then the set $$\{p(\iota^2) \iota : p \hbox{ a polynomial with complex coefficients}\}$$ is dense in $C_0 (\Omega)$. Therefore, given $\varepsilon>0$ there exists a polynomial $p$ with complex coefficients such that for $q(\lambda) = p(\lambda^2) \lambda$ we have $\sup_{t\in \Omega} |q(t) - f(t) |<\varepsilon$. By the continuity of $q$ --via the continuity of the triple product-- we can find $\delta > 0$ such that $\| q_t(a)- q_t(b)\| < \varepsilon$ for all $a,b\in E$ with $\|a-b\|<\delta$ and $\|a\|,$ $\|b\|\leq \max \Omega$. Taking $a,b\in \mathcal{E}_{\Omega}$ with $\|a-b\|<\delta$ we have $$\|f_t(a) -f_t(b) \|\leq \|f_t(a) -q_t(a) \| + \|q_t(a) -q_t(b) \| + \|q_t(b) -f_t(b) \|< 3\varepsilon.$$ The last statement is clear.
\end{proof}

The next corollary is our tool to avoid polar decompositions in the Jordan setting.

\begin{corollary}\label{c continuity of the range tripotent} Let $M$ be a unital JB$^*$-algebra. Then the mapping $$a\in M^{-1}\mapsto r(a)\in \mathcal{U} (M) \hbox{ is continuous.}$$
\end{corollary}

\begin{proof} It is not hard to check, from the properties of the continuous triple functional calculus, that for each von Neumann regular element $a\in E$ the identity $ r(a) = Q(a^{\dag}) (a^{[2]})$ holds.\smallskip

We recall that the set of invertible elements in a unital Jordan-Banach algebra $M$ is an open set, and the mapping $a\in M^{-1}\mapsto a^{-1}\in M$ is continuous (cf. \cite[Theorem 4.1.7 and Proposition 4.1.6]{Cabrera-Rodriguez-vol1}). Therefore, if $M$ is a unital JB$^*$-algebra, the mapping $a\in M^{-1}\mapsto a^{\dag} = (a^{-1})^*\in M^{-1}$ is continuous. Consequently, by the continuity of the triple product, the mapping $$a\in M^{-1} \mapsto r(a) = Q(a^{\dag}) (a^{[2]})$$ is continuous. Let us finally note that $r(a)\in \mathcal{U} (M)$ for all $a\in M^{-1}$ (cf. \eqref{eq range of invertible is unitary}).
\end{proof}

We have already employed the notion of $u$-isotope associated with each unitary element $u$ in a unital JB$^*$-algebra $M$. This device also makes sense for invertible elements in $M$. More concretely, for each invertible element $c$ in a unital Jordan--Banach algebra $M$ then the vector space $M$ becomes a Jordan--Banach algebra $M_{(c)}$ with unit element $c^{-1}$ and quadratic operators \begin{equation}\label{eq quadratic operation on isotopes} U^{(c)}_a = U_a U_c.
 \end{equation} The linear Jordan product in $M_{(c)}$ is actually given by $x \circ_{c} y = x\circ (c\circ y) + y\circ (c \circ x) - c\circ (x\circ y)$ (cf. \cite[\S 2, 1.7]{Jacobson81}). We should note that, according to the notation introduced in Remark \ref{r isotopes remark}, $M(w)= M_{(w^*)}$.\smallskip

As O. Loos observed in \cite[page 112]{Loos96}, ``\emph{... passing to an isotope is in some sense a substitute for the left multiplication with an invertible element in an associative Banach algebra, an operation which has no direct Jordan analogue. It is therefore important to find isotopy--invariant
objects in a Jordan--Banach algebra.}''\smallskip


The algebraic characterization of the principal component of the set of unitaries in a unital JB$^*$-algebra reads as follows.

\begin{theorem}\label{t principal component as product of exp} Let $M$ be a unital JB$^*$-algebra, let $\mathcal{U}^0(M)$ denote the principal component of the set of unitaries in $M$ and let $u$ be an element in $\mathcal{U}(M)$. Then the following statements are equivalent:
\begin{enumerate}[$(a)$]\item $u\in M^{-1}_{\textbf{1}}\cap\mathcal{U} (M)$;
\item There exists a continuous path $\Gamma : [0,1]\to \mathcal{U} (M)$ with $\Gamma (0) = \textbf{1}$ and $\Gamma(1) = u$;
\item $u\in \mathcal{U}^0(M)$;
\item $u = U_{e^{i h_n}} \cdots U_{e^{i h_1}}(\textbf{1}),$ for some $n\in \mathbb{N}$, $h_1,\ldots,h_n\in M_{sa}$;
\item There exists $w\in \mathcal{U}^0(M)$ such that $\|u-w\|<2$.
\end{enumerate}
Consequently,
\begin{equation}\label{eq algebraic charact principal component unitaries JBstar}\begin{aligned} \mathcal{U}^0(M) &=  M^{-1}_{\textbf{1}}\cap\mathcal{U} (M)\\
&=\left\lbrace  U_{e^{i h_n}}\cdots U_{e^{i h_1}}(\textbf{1}) \colon n\in \mathbb{N}, \ h_j\in M_{sa} \ \forall\ 1\leq j \leq n  \right\rbrace \\
&= \left\lbrace  u\in \mathcal{U} (M) : \hbox{ there exists } w\in \mathcal{U}^0(M) \hbox{ with } \|u-w\|<2 \right\rbrace
\end{aligned}\end{equation} is analytically arcwise connected.
\end{theorem}

\begin{proof} $(a)\Rightarrow (b)$ Suppose $u\in M^{-1}_{\textbf{1}}\cap\mathcal{U} (M)$. Since by Loos' Theorem  \ref{t Loos principal component invertible}, the set $M^{-1}_{\textbf{1}}$ is analytically arcwise connected, there exists a continuous path $\gamma : [0,1]\to M^{-1}_{\textbf{1}}$ satisfying $\gamma(0) = \textbf{1}$ and $\gamma(1) = u$. Corollary \ref{c continuity of the range tripotent} assures that the mapping $\Gamma : [0,1]\to \mathcal{U} (M),$ $\Gamma (t) := r(\gamma(t))$ is continuous with $\Gamma (0) = \textbf{1}$ and $\Gamma(1) = u$. \smallskip

The implication $(b)\Rightarrow (c)$ is clear.\smallskip

$(c)\Rightarrow (a)$ It follows from $(a)\Rightarrow (b)\Rightarrow (c)$ that $M^{-1}_{\textbf{1}}\cap\mathcal{U} (M)\subseteq \mathcal{U}^0(M)$. We know from Loos' Theorem  \ref{t Loos principal component invertible} that $M^{-1}_{\textbf{1}}$ is a clopen subset of $M^{-1}$, and thus $M^{-1}_{\textbf{1}}\cap\mathcal{U} (M)$ is a clopen subset of $\mathcal{U} (M)$ which contains the unit and is contained in $\mathcal{U}^0(M)$. Having in mind that $\mathcal{U}^0(M)$ is a connected set, we deduce that $M^{-1}_{\textbf{1}}\cap\mathcal{U} (M) = \mathcal{U}^0(M).$\smallskip

We have proved that statements $(a)$, $(b)$ and $(c)$ are equivalent and $\mathcal{U}^0(M)= M^{-1}_{\textbf{1}}\cap\mathcal{U} (M)$ is analytically arcwise connected.\smallskip

$(b)\Rightarrow (d)$  Suppose there exists a continuous path $\Gamma : [0,1]\to \mathcal{U} (M)$ with $\Gamma (0) = \textbf{1}$ and $\Gamma(1) = u$. Let us find, by continuity and compactness, $0=t_0 < t_1 < \ldots < t_n < 1= t_{n+1}$ such that $\|\Gamma (t_i) - \Gamma (t_{i+1})\|<2$ for all $0\leq i\leq n$.  Set $u_i=\Gamma(t_i)$ for all $0\leq i \leq n+1$. By applying Lemma \ref{l untaries at distance <2 in a unital JB$^*$-algebra} to the element $u_1$ --which satisfies $\|u_1 -\textbf{1}\|<2$-- we deduce the existence of $h_1\in M_{sa}$ such that $U_{e^{i {h_{1}}}} (\textbf{1}) = u_1$.  Since $2> \|u_1 - u_2\| = \|\textbf{1} - U_{ e^{-i {h_{1}}}} ( u_2)\|$, a new application of Lemma \ref{l untaries at distance <2 in a unital JB$^*$-algebra} shows the existence of $h_{2}\in M_{sa}$ satisfying $  U_{ e^{-i {h_{1}}}} ( u_2) =  U_{ e^{i {h_2}}} ( \textbf{1}), $ or equivalently, $ u_2  = U_{ e^{i {h_{1}}}}  U_{ e^{i {h_2}}} ( \textbf{1})$. Suppose, by an induction argument on $n$, that we have found $h_1,\ldots,h_n\in M_{sa}$ with $u_{n} = U_{e^{i h_1}} \cdots U_{e^{i h_n}}(\textbf{1}).$ As before, the condition $$2>\| u_{n+1} - u_{n} \| =  \| U_{e^{-i h_n}} \cdots  U_{e^{-i h_1}} ( u_{n+1} ) - \textbf{1} \|,$$ implies, via Lemma \ref{l untaries at distance <2 in a unital JB$^*$-algebra}, the existence of $h_{n+1}\in M_{sa}$ such that $u = u_{n+1} =  U_{e^{i h_{1}}}  U_{e^{i h_2}} \cdots U_{e^{i h_{n+1}}}(\textbf{1}).$\smallskip

The implication $(d)\Rightarrow (c)$ is also easy because given $u = U_{e^{i h_n}} \cdots U_{e^{i h_1}}(\textbf{1}),$ with $n\in \mathbb{N}$, $h_1,\ldots,h_n\in M_{sa},$ the mapping $\Gamma: [0,1]\to \mathcal{U} (M) \subseteq M^{-1}$, $\Gamma (t) = U_{e^{i t h_n}} \cdots U_{e^{i t h_1}}(\textbf{1})$ is an analytic curve with $\Gamma (0) = \textbf{1}$ and $\Gamma (1) = u$. This also proves the second equality in \eqref{eq algebraic charact principal component unitaries JBstar}.\smallskip

The implication $(c)\Rightarrow (e)$ is clear.\smallskip

Suppose finally that statement $(e)$ holds, that is, there exists $w\in \mathcal{U}^{0} (M)$ satisfying $\|w-u\|<2$. By the first two equalities in \eqref{eq algebraic charact principal component unitaries JBstar}, $w = U_{e^{i h_n}} \cdots U_{e^{i h_1}}(\textbf{1})$ for some $h_1, \ldots, h_{n}$ in $M_{sa}$. A similar argument to that in the proof of $(b)\Rightarrow (d)$ implies that $u = U_{e^{i h_{n+1}}} U_{e^{i h_n}} \cdots U_{e^{i h_1}}(\textbf{1})\in \mathcal{U}^{0} (M)$.
\end{proof}

In case that $u$ is a unitary in an associative unital C$^*$-algebra $A$, the left multiplication operator $L_u: A\to A$ is a surjective linear isometry mapping the unit to $u$. Consequently, $L_u (\mathcal{U}^0(A)) = u  \mathcal{U}^0(A)$ is precisely the connected component of $\mathcal{U}(A)$ containing $u$. However, the left multiplication operator admits no direct Jordan analogue. We shall see next that assuming that a unitary element $u$ in a unital JB$^*$-algebra $M$ writes as the square of another unitary, the connected component containing it can be described.

\begin{remark}\label{r unitaries which are squares} Let $u$ be a unitary element in a unital JB$^*$-algebra $M$. We cannot always guarantee that $u = e^{i h}$ for some $h\in M_{sa}$ {\rm(}cf. \cite[Exercise 4.6.9]{KR1}{\rm)}. Suppose that $u$ satisfies the weaker condition that $u = v^2$ for some unitary $v\in M$ --this automatically holds if $u = e^{i h}$. Let $\mathcal{U}^c(M)$ denote the connected component of $\mathcal{U}(M)$ containing $u$. Since the mapping $U_v: M\to M(u)$ is a unital surjective linear isometry --and hence a Jordan $^*$-isomorphism between these two algebras-- $U_{v} (\mathcal{U}^{0}(M))$ is a connected component of $\mathcal{U} (M) =\mathcal{U} (M(u))$ and contains $u$. It then follows from Theorem \ref{t principal component as product of exp} that $$ \mathcal{U}^c (M) = U_{v} (\mathcal{U}^{0}(M)) = \left\lbrace U_{v} U_{e^{i h_n}}\cdots U_{e^{i h_1}}(\textbf{1}) \colon n\in \mathbb{N}, \ h_j\in M_{sa} \ \forall\ 1\leq j \leq n  \right\rbrace.$$

The existence of a unitary square root does not necessarily hold for any unitary element, even in the setting of commutative C$^*$-algebras. Let $\mathbb{T}$ denote the unit sphere of the complex plane. The principal component of $\mathcal{U}(C(\mathbb{T}))$ is precisely the subgroup ${\rm\exp}(i C_{\mathbb{R}}(\mathbb{T}))$, which is the set of all functions $u: \mathbb{T}\to \mathbb{T}$ which are deformable or homotopic to the unit element {\rm(}cf. \cite[Exercise 4.6.7]{KR1}{\rm)}. The quotient group $\mathcal{U}(C(\mathbb{T}))/ \mathcal{U}^0(C(\mathbb{T})) = \mathcal{U}(C(\mathbb{T}))/ \exp(i C_{\mathbb{R}}(\mathbb{T}))$ --known as the Bruschlinsky group-- identifies with the ring $\mathbb{Z}$ of integers {\rm(}see \cite[\S II.3]{Hu59} or \cite[Exercise I.11.3]{Tak}{\rm)}. It is known that the maps $u_1, u_2:\mathbb{T}\to \mathbb{T},$ $u_1(\lambda)= \lambda$ and $u_2(\lambda)= \lambda^2$ are not in the principal component of $\mathcal{U}(C(\mathbb{T}))$, they are actually in two different connected components {\rm(}cf. \cite[\S II, Lemma 3.2]{Hu59}{\rm)}. Therefore, $u_2\notin \exp(i C_{\mathbb{R}}(\mathbb{T}))$ but $u_2 = u_1^2$ admits a unitary square root.
\end{remark}

The hypothesis assuring that a unitary element $u$ in a JB$^*$-algebra $M$ admits a square root has been already considered in the literature, for example, in \cite[Problem 5.1 and Lemma 5.2]{BraKaUp78}, \cite{HatMi2000}, \cite{ChiKasKawValov2005} and \cite{MatNagYam2008}. The good behavior exhibited by unital C$^*$-algebras, where the connected components of the unitary group can be derived from the principal component by just multiplying by a determined unitary on the left or on the right, fails in the setting of unital JB$^*$-algebras, not only due to the lacking of a left (or right) multiplication operator, but for the possibility that two different connected components of the unitary set of a unital JB$^*$-algebra are not isometric as metric spaces (see Remark \ref{r existence of unitaries with non isomorphic nor isometric connected components in the Jordan setting}).\smallskip

The lacking of a left multiplication operator can be also overcome by passing to an isotope.

\begin{remark}\label{r connected components which are not the pricipal one} Let $w$ be a unitary element in a unital JB$^*$-algebra $M$. Suppose $\mathcal{U}^c(M)$ is a connected component of $\mathcal{U} (M)$, $w\in \mathcal{U}^c(M)$ and $M^{-1}_{w}$ is the connected component of $M^{-1}$ containing $w$. Let $M(w)$ denote the $w$-isotope admitting $w$ as unit element. Since $\mathcal{U}(M) = \mathcal{U}(M(w))$ and $M^{-1} = (M(w))^{-1}$ as sets {\rm(}cf. Remark \ref{r isotopes remark}{\rm)}, $\mathcal{U}^c(M)$ and $M^{-1}_{w}$ are the principal components of $\mathcal{U}(M(w))$ and $(M(w))^{-1}$, respectively. Therefore the previous Theorem \ref{t principal component as product of exp} applied to the $w$-isotope $M(w)$ gives \begin{equation}\label{eq algebraic charact secondary component unitaries JBstar} \begin{aligned} \mathcal{U}^c (M)&=  M^{-1}_{w}\cap\mathcal{U} (M)\\
&=\left\lbrace  U^{(w)}_{e_w^{i h_n}}\cdots U^{(w)}_{e_w^{i h_1}}(w) \colon n\in \mathbb{N}, \ h_j\in (M(w))_{sa}, \ \forall\ 1\leq j \leq n  \right\rbrace \\
&= \left\lbrace  u\in \mathcal{U} (M) : \hbox{ there exists } v\in \mathcal{U}^c(M) \hbox{ with } \|u-v\|<2 \right\rbrace
\end{aligned}
\end{equation} is analytically arcwise connected. We can actually take any $\widetilde{w}\in \mathcal{U}^c(M)$ and the same conclusion holds for the connected component $\mathcal{U}^c(M)$ and the $\widetilde{w}$-isotope  $M(\widetilde{w})$.
\end{remark}

It will be useful to establish several further consequences of the previous Theorem \ref{t principal component as product of exp}.

According to \cite{Pe2020}, given a unital JB$^*$-algebra $M$, a subset $\mathfrak{M} \subseteq M^{-1}$ will be called a \emph{quadratic subset} if $U_{\mathfrak{M}} (\mathfrak{M}) \subseteq \mathfrak{M}$. Let us observe that $\mathcal{U} (M)$ is a self-adjoint subset of $M$.\smallskip

\begin{proposition}\label{p self adjoint quadratic subset princ comp} Let $M$ be a unital JB$^*$-algebra. Then $\mathcal{U}^{0} (M)$ is a self-adjoint quadratic subset of $\mathcal{U} (M)$, that is, for all $u,w\in \mathcal{U}^0 (M)$ the elements $w^*$ and $U_w (u)$ belong to $\mathcal{U}^0 (M)$. Furthermore, $\mathcal{U}^{0} (M)$ is the smallest quadratic subset of $\mathcal{U} (M)$ containing the set $e^{i M_{sa}}$. In particular, $$U_{\mathcal{U}^{0} (M)} \left( \mathcal{U}^{0} (M) \right)\! = \mathcal{U}^{0} (M) =\left\lbrace  U_{e^{i h_n}}\cdots U_{e^{i h_1}}(e^{i h_0}) \colon n\in \mathbb{N}, \ h_j\in M_{sa} \ \forall\ 0\leq j \leq n  \right\rbrace.$$
\end{proposition}

\begin{proof} Since $u,w\in \mathcal{U}^0(M),$ Theorem \ref{t principal component as product of exp} guarantees the existence of hermitian elements $h_m,\dots,h_1, k_n, \ldots, k_1\in M_{sa}$ such that $w=U_{e^{i h_m}}\cdots U_{e^{i h_1}}(\textbf{1})$ and $u=U_{e^{i k_n}}\cdots U_{e^{i k_1}}(\textbf{1})$ with $n,m\in \mathbb{N}$. Concerning the involution we have
$$ w^*= (U_{e^{i h_m}}\cdots U_{e^{i h_1}}(\textbf{1}))^* =\hbox{(by \eqref{eq U-op adjoint})} = U_{e^{-i h_m}}\cdots U_{e^{-i h_1}}(\textbf{1}),$$  which proves that $w^*\in \mathcal{U}^0(M)$ (cf. Theorem \ref{t principal component as product of exp}).\smallskip

Concerning the remaining property, we apply the identity in \eqref{eq generalised fundamental identity UaUbUa} to deduce that
$$\begin{aligned}
U_w (u)&=U_w\left( U_{e^{i k_n}}\cdots U_{e^{i k_1}}(\textbf{1}) \right) = U_{U_{e^{i h_m}}\cdots U_{e^{i h_1}}(\textbf{1})}\left( U_{e^{i k_n}}\cdots U_{e^{i k_1}}(\textbf{1}) \right) \\
&= U_{e^{i h_m}}\cdots U_{e^{i h_1}} U_{\textbf{1}} U_{e^{i h_1}}\cdots U_{e^{i h_m}}\left( U_{e^{i k_n}}\cdots U_{e^{i k_1}}(\textbf{1}) \right)\\
&=U_{e^{i h_m}}\cdots U_{e^{i h_1}} U_{e^{i h_1}}\cdots U_{e^{i h_m}} U_{e^{i k_n}}\cdots U_{e^{i k_1}}( \textbf{1}),
\end{aligned}$$ which in view of Theorem \ref{t principal component as product of exp} assures that $U_w (u)\in \mathcal{U}^0(M)$.\smallskip

Fix an arbitrary $w$ in $\mathcal{U}^0(M)$. Since $U_w$ is an invertible mapping with inverse $U_{w^*},$ the above properties show that $w^*\in \mathcal{U}^0(M)$ and
$$\mathcal{U}^0(M) = U_w  U_{w^*} (\mathcal{U}^0(M)) \subseteq  U_w (\mathcal{U}^0(M)) \subseteq \mathcal{U}^0(M).$$

Finally, suppose that $\mathcal{V}$ is a quadratic subset of $\mathcal{U} (M)$ containing the set $e^{i M_{sa}}$. Clearly for each $h_1,\ldots, h_n$ in $M_{sa}$ an induction argument shows that $$U_{e^{i h_n}} \ldots U_{e^{i h_1}} (\textbf{1}) \in \mathcal{V},$$ and thus, by Theorem \ref{t principal component as product of exp}, $\mathcal{U}^0(M)\subseteq \mathcal{V}$. The rest is clear.
\end{proof}

We conclude this section with a remark and an open problem.

\begin{remark}\label{r connected component building from u as an pen question} Let $M$ be a unital JB$^*$-algebra and let $u$ be an element in $\mathcal{U} (M)\backslash \mathcal{U}^0(M)$. Suppose $\mathcal{U}^c (M)$ is the connected component of $\mathcal{U} (M)$ containing $u$. The set $$\mathcal{E}_u := \left\{ U_{e^{i h_1}}\ldots U_{e^{i h_n}} (u) : \begin{array}{c}
                                                                n\in \mathbb{N}, \ h_j\in M_{sa} \\
                                                                \hbox{ for all } 1\leq j \leq n
                                                              \end{array}
\right\}$$ is arcwise connected and contains $u$. Thus, $\mathcal{E}_u\subseteq \mathcal{U}^c (M)$. We do not know if the reverse inclusion holds, in general. If $M= A$ is a unital C$^*$-algebra both inclusions hold, a fact which is essentially due to the existence of left and right multiplicaiton operators. Namely, since $\mathcal{U}^0(A) = \{ e^{i h_1} \cdot \ldots \cdot e^{i h_n} : n\in \mathbb{N}, \ h_1, \ldots, h_n\in A_{sa} \}$ {\rm(}cf. \eqref{eq algebraic descirption of the principal component for unital Cstar algebras}{\rm)}, and the right multiplication operator $R_u: A\to A$, $x\mapsto x u$ is a surjective linear isometry mapping the unit to $u$, we deduce that $$\left\{ e^{i h_1} \cdot \ldots \cdot e^{i h_n} u : n\in \mathbb{N}, \ h_1, \ldots, h_n\in A_{sa} \right\}= R_u (\mathcal{U}^0(A)) = \mathcal{U}^c(A).$$

Similarly, the products of the form $w e^{i h_n} \cdot \ldots \cdot e^{i h_1}$ lie in $\mathcal{U}^c (A)$ for all $w\in \mathcal{U}^c(A)$, $h_1, \ldots, h_n\in A_{sa}$. So, given $h_1, \ldots, h_n\in A_{sa}$, the mapping $R_{e^{i h_1} \cdot \ldots \cdot e^{i h_n}}$ is a surjective linear isometry on $A$ mapping $\mathcal{U}^c(A)$ onto itself, and hence $$\begin{aligned}
\mathcal{U}^c(A)  &= \left\{ e^{i h_1} \cdot \ldots \cdot e^{i h_n} u e^{i h_n} \cdot \ldots \cdot e^{i h_1} : n\in \mathbb{N}, \ h_1, \ldots, h_n\in A_{sa} \right\} \\
& = \left\{ U_{e^{i h_1}}\ldots U_{e^{i h_n}} (u) : \begin{array}{c}
                                                                n\in \mathbb{N}, \ h_j\in A_{sa} \\
                                                                \hbox{ for all } 1\leq j \leq n
                                                              \end{array}
\right\}.
\end{aligned} $$\smallskip

Despite we do not know if the equality $\mathcal{E}_u= \mathcal{U}^c (M)$ holds in general, the inclusion $\mathcal{E}_u\subseteq \mathcal{U}^c (M)$ is enough to prove the following rule concerning the $U$-operators \begin{equation}\label{eq inner main other main} U_{\mathcal{U}^0 (M)} \left( \mathcal{U}^c (M)\right) = \mathcal{U}^c (M), \hbox{ for all connected component } \mathcal{U}^c (M) \hbox{ of } \mathcal{U} (M).
\end{equation} Indeed, by Theorem \ref{t principal component as product of exp}, every $w\in \mathcal{U}^{0}(M)$ is of the form $w = U_{e^{i h_n}}\cdots U_{e^{i h_1}}(\textbf{1})$ with $n\in \mathbb{N}$, $h_j\in M_{sa}$. Therefore, given $u\in \mathcal{U}^c (M)$, it follows from \eqref{eq generalised fundamental identity UaUbUa} that $$ U_w (u) = U_{U_{e^{i h_n}}\cdots U_{e^{i h_1}}(\textbf{1})} (u) = U_{e^{i h_n}} \ldots U_{e^{i h_1}} U_{e^{i h_1}} \ldots U_{e^{i h_n}} (u)\in \mathcal{E}_u\subseteq \mathcal{U}^c (M).$$\smallskip

We finally observe that the equality $$\bigcup_{w\in \mathcal{U}^{c}(M)} \mathcal{E}_w = \mathcal{U}^{c}(M)$$ follows trivially from the above arguments.
\end{remark}

\section{Surjective isometries between principal components}\label{sec: surj isom}

According to the notation in \cite{HatHirMiuMol2012,CuPe20}, let $\mathcal{G}$ be a group and let $(X,d)$ be a non-trivial metric space such that $X$ is a subset of $\mathcal{G}$ and \begin{equation}\label{eq condition p HatoriMiuraMolnar}\hbox{ $y x^{-1} y\in X$ for all $x,y\in X$}
\end{equation} (note that we are not assuming that $X$ is a subgroup of $\mathcal{G}$). We need to recall some basic definitions.

\begin{definition}\label{def condition B} Let us fix $a,b$ in $X$. We shall say that condition $B(a,b)$ holds for $(X,d)$ if the following properties are satisfied:\begin{enumerate}[$(B.1)$]\item For all $x,y\in X$ we have $d(b x^{-1} b, b y^{-1} b) = d(x,y).$
		\item There exists a constant $K>1$ satisfying $$d(b x^{-1} b,x)\geq K d(x,b),$$ for all $x\in L_{a,b} =\{x\in X : d(a,x) = d(b a^{-1} b, x) = d(a,b)\}.$
	\end{enumerate}
\end{definition}

\begin{definition}\label{def condition C2} Let us fix $a,b\in X$. We shall say that condition $C_1(a,b)$ holds for $(X,d)$ if the following properties are satisfied:\begin{enumerate}[$(C.1)$]\item For every $x\in X$ we have $a x^{-1} b, b x^{-1} a\in X$;
		\item $d(a x^{-1} b, a y^{-1} b) = d(x,y)$, for all $x,y\in X$.
	\end{enumerate}\smallskip
	
	\noindent We shall say that condition $C_2(a,b)$ holds for $(X,d)$ if there exists $c\in X$ such that $c a^{-1} c =b$ and $d(c x^{-1} c, c y^{-1} c) = d(x,y)$ for all $x,y\in X$.
\end{definition}

An element $x\in X$ is called \emph{2-divisible} if there exists $y\in X$ such that $y^2 =x$. $X$ is called \emph{2-divisible} if every element in $X$ is {2-divisible}. Furthermore, $X$ is called \emph{2-torsion free} if it contains the unit of $\mathcal{G}$ and the condition $x^2 =1$ with $x\in X$ implies $x =1$.\smallskip

As in previous references (cf. \cite{HatMol2012,HatMol2014,CuPe20b}) our arguments here will rest on the following technical result, which has been borrowed from \cite{HatHirMiuMol2012}. 

\begin{theorem}\label{t HatHirMiuMol Thm 24 for metric spaces}\cite[Theorem 2.4]{HatHirMiuMol2012} Let $(X, d_X)$ and $(Y, d_Y )$ be two metric spaces. Pick two points $a,c \in X$. Suppose that $\varphi : X \to X$ is a distance preserving map such that $\varphi(c) = c$ and $\varphi\circ \varphi$ is the identity map on $X$. Let $$L_{a,c} = \{x \in X : d_X (a, x) = d_X(\varphi(a), x) = d_X(a, c)\}.$$
	Suppose that there exists a constant $K > 1$ such that $d_X (\varphi(x), x) \geq K d_X(x, c)$ holds for every $x \in L_{a,c}$. If $\Delta$ is a bijective distance preserving map from $X$ onto $Y$, and $\psi$ is a bijective distance preserving map from $Y$ onto itself such that $\psi(\Delta (a)) = \Delta (\varphi(a))$ and $\psi(\Delta (\varphi(a))) = \Delta (a),$ then we have $\psi(\Delta (c)) =\Delta (c).$ 
\end{theorem}

Theorem 2.9 in \cite{CuPe20b} decodes the main properties of a surjective isometry between the sets of unitary elements in two unital JB$^*$-algebras. In the next theorem we establish an analogous version of the just quoted result but for surjective isometries between the principal components of the sets of unitary elements. The original proof requires a slightly modified argument.

\begin{theorem}\label{t Delta preserves inverted triple products in U0M} Let $\Delta: \mathcal{U}^0 (M)\to \mathcal{U}^0 (N)$ be a surjective isometry, where $M$ and $N$ are unital JB$^*$-algebras. Suppose $u,v\in \mathcal{U}^0 (M)$ with $\|u-v\|<\frac12$. Then the following statements are true:
\begin{enumerate}[$(1)$]\item The Jordan version of condition $B(u,v)$ holds for $\mathcal{U}^0 (M)$, that is, the following statements are satisfied:\smallskip
\begin{enumerate}[$(a)$]
			\item For all $x,y\in \mathcal{U}^0 (M)$ we have $$\left\|U_{v}(x^{-1})- U_{v} (y^{-1})\right\|=\left\|U_{v}(x^{*})- U_{v} (y^{*})\right\|= \|x^*-y^*\| = \|x-y\|.$$
			\item The constant $K= 2 - 2 \|u-v\|>1$ satisfies that $$ \left\| U_v ( w^*)- w \right\|=\left\| U_v ( w^{-1})- w \right\|\geq K \| w- v\|,$$ for all $w$ in the set $$L^0_{u,v} =\{w\in \mathcal{U}^0 (M)  : \|u-w\| = \| U_{v} (u^{-1})- w\|= \| U_{v} (u^*)- w\| = \|u-v\|\}.$$
		\end{enumerate}\smallskip
		\item The Jordan version of condition $C_2 (\Delta(u),\Delta(U_v (u^{*})))$ holds for $\mathcal{U}^0 (N)$, that is,
		there exists $w\in \mathcal{U}^0 (N)$ satisfying \begin{equation}U_w (\Delta(u)^*) = \Delta\left( U_v (u^*) \right),\nonumber
		\end{equation}
		and
		$$\| U_w (x^{*}) -U_w (y^{*})\|  = \|x-y\|, \quad \forall x,y\in \mathcal{U}^0 (N).$$\smallskip
		\item The equalities $\Delta(U_{v} (u^{*})) =\Delta(U_{v} (u^{-1})) =U_{\Delta(v)} (\Delta(u)^{-1}) =U_{\Delta(v)} (\Delta(u)^*)$ hold.
	\end{enumerate}
\end{theorem}

\begin{proof} Let us note that Theorem 2.9 in \cite{CuPe20b} proves that the same conclusion hold for every surjective isometry $\Delta: \mathcal{U}^0 (M)\to \mathcal{U}^0 (N)$, with the unique exception that $\mathcal{U}^0 (M)$ is replaced with $\mathcal{U} (M)$.\smallskip

Since $\mathcal{U}^0(M)\subseteq \mathcal{U}(M)$, statement $(1)$ follows directly from \cite[Theorem 2.9{\rm$(a)$}]{CuPe20b} (see also \cite[Lemma 2.8]{CuPe20b}).\smallskip
	
	$(2)$ Let us consider the $u$-isotope $M(u)$, and observe that $v\in\mathcal{U} (M(u))$ (cf. Remark \ref{r isotopes remark}$(b)$).  Let $\mathcal{B}$ denote the JB$^*$-subalgebra of $M(u)$ generated by $v$ and the unit of $M(u)$, that is, $u$. The Shirshov-Cohn theorem assures that $\mathcal{B}$ is a JB$^*$-subalgebra of $B(H)$, for some complex Hilbert space $H$, and we can further assume that $u$ is the unit in $B(H)$ (cf. \cite[Theorem 7.2.5]{HOS} and \cite[Corollary 2.2]{Wri77}). Moreover, $u$ can be identified with the unit of $B(H)$. We shall use juxtaposition for the product in $B(H)$, and the symbol $\sharp$ for its involution.\smallskip
	
It is worth noting that when $M$ and $\mathcal{B}$ are regarded as JB$^*$-triples, the triple product in $\mathcal{B}$ is precisely  $\J{\cdot}{\cdot}{\cdot}_M|_{\mathcal{B}\times \mathcal{B}\times \mathcal{B}}$. At the same time, the triple product in $\mathcal{B}$ can be expressed in terms of the associative product of $B(H)$. By Proposition \ref{p self adjoint quadratic subset princ comp} the element $U_v (u^*)$ lies in $\mathcal{U}^0(M)$, and hence we can apply $\Delta$ at this element and compute the next distances:
	$$\begin{aligned} \| \Delta(u) - \Delta(U_v (u^*))\|_N &=  \| u - U_v (u^*)\|_M =   \| u - \{ v, u,v\}_M\|_M = \| u - \{ v, u,v\}_{\mathcal{B}}\|_{\mathcal{B}}  \\
	&= \| u - \{ v, u,v\}_{_{B(H)}}\|_{_{B(H)}} =  \| u - v u^{\sharp} v\|_{_{B(H)}} \\
	&=  \| u - v v\|_{_{B(H)}} =  \| v^{\sharp} - v\|_{_{B(H)}} \leq   \| v^{\sharp} - u \|_{_{B(H)}} + \| u- v\|_{_{B(H)}} \\
	&= 2 \| u - v\|_{_{B(H)}}= 2 \| u - v\|_{\mathcal{B}} =  2 \| u - v\|_M < 1.
	\end{aligned}.$$ 

Actually, there is another argument which avoids any use of the Shirshov-Cohn theorem. Namely, since $\mathcal{B}$ is the JB$^*$-subalgebra of the $u$-isotope $M(u)$ generated by the unitary $v$ --which is a unitary in $M(u)$, see Remark \ref{r isotopes remark}-- and the unit of $M(u)$ --i.e. $u$-- it must be isometrically Jordan $^*$-isomorphic to a unital commutative C$^*$-algebra, that is, to some $C(\Omega)$ for an appropriate compact Hausdorff space $\Omega,$ and under this identification, $u$ corresponds to the unit  (cf. \cite[3.2.4. The spectral theorem]{HOS} or \cite[Proposition 3.4.2 and Theorem 4.1.3$(v)$]{Cabrera-Rodriguez-vol1}). So, working in this commutative unital C$^*$-algebra $\mathcal{B}\cong C(\Omega)$, we can compute the above distances in an even easier setting.\smallskip
	
We are in a position to apply Lemma 2.2$(b)$ in \cite{CuPe20b} to assure the existence of a unitary $w\in \mathcal{U} (N)$ satisfying $$U_w ( \Delta(u)^* ) = \Delta\left( U_v (u^*) \right).$$ However, by hypotheses $\| \Delta(u) - \Delta(v)\| =\|u-v\|<\frac12$, and consequently, \cite[Lemma 2.2$(b)$]{CuPe20b} also guarantees that $$\|w- \Delta(v) \|\leq \|w- \Delta(u) \|+\|\Delta(u)-\Delta(v) \|\leq \| \Delta(u) - \Delta(U_v (u^*))\| +\frac12  < 2.$$ Since $\Delta(v)\in \mathcal{U}^0 (N)$, Theorem \ref{t principal component as product of exp} implies that $w$ belongs to the principal component of $\mathcal{U}(N)$, that is, $w\in \mathcal{U}^0 (N)$.\smallskip

Now, it follows from Remark \ref{r isotopes remark}$(c)$ that
	$$\| U_w (x^{*}) -U_w (y^{*})\|=\| U_w (x^{-1}) -U_w (y^{-1})\| = \|x^{-1}- y^{-1}\| = \|x^* - y^* \| = \|x-y\|,$$ for all $x,y\in \mathcal{U}^0 (N)\subseteq \mathcal{U}(N)$. We can hence conclude that the Jordan version of $C_2 (\Delta(u),\Delta(U_u (v^{*})))$ holds for $\mathcal{U}^0 (N)$.\smallskip


$(3)$ Let $w\in \mathcal{U}^0 (N)$ be the element given by $(2)$.\smallskip

Let us define the mappings $\varphi: \mathcal{U}^0 (M)\to \mathcal{U}^0 (M)$ and $\psi: \mathcal{U}^0 (N)\to \mathcal{U}^0 (N)$ given by
$$\varphi (x) := U_{v} (x^{-1}) = U_{v} (x^{*}), \ \ x\in \mathcal{U}^0(M),$$ and
$$\psi (y) := U_{w} (y^{-1}) = U_{w} (y^{*}), \ \ y\in \mathcal{U}^0(N),$$ respectively. They are well defined thanks to Proposition \ref{p self adjoint quadratic subset princ comp}.\smallskip

We aim to apply Theorem \ref{t HatHirMiuMol Thm 24 for metric spaces} to $\varphi, \psi$ and $\Delta$ at $a= u$ and $c = v$. For this purpose let us observe that, since $v$ and $w$ are unitaries lying in the respective principal components, $\varphi$ and $\psi$ are distance preserving bijections (cf. Remark \ref{r isotopes remark}$(c)$).\smallskip
	
On the other hand, the element $v$ is invertible with $v^*$ as its inverse in $M$. That implies that $\varphi (v)=U_{v} (v^{-1})=U_{v} (v^*) = v$ and $$\varphi\circ \varphi(x)=U_{v} ((U_{v} (x^{-1}))^{-1})=U_{v} ((U_{v} (x^*))^*)= U_{v} (U_{v^*} (x))=U_{v} (U_{v^{-1}} (x))= x,$$ for every $x\in \mathcal{U}^0(M)$ (cf. \eqref{eq simple U-op adjoint}). That is, $\varphi\circ \varphi$ is the identity mapping on $\mathcal{U}^0 (M)$. By applying that $U_w ( \Delta(u)^* ) = \Delta\left( U_v (u^*) \right)$ we have $$\begin{aligned}
	\psi (\Delta(u)) &= U_w (\Delta(u)^*) = \Delta\left( U_v (u^*) \right) =\Delta\left( \varphi (u) \right), \hbox{ and }\\
	\psi \left(\Delta\left( \varphi (u) \right)\right) &= U_{w} \left((\Delta \left( U_{v} (u^*)\right))^* \right)= \left( U_{w^*} \left(\Delta \left( U_{v} (u^*)\right)\right)\right) ^* = \Delta(u)^{**} = \Delta (u).
	\end{aligned}$$
	
Therefore the mappings  $\varphi, \psi,$ $\Delta,$ $a= u$ and $c=v$ satisfy all hypotheses in Theorem \ref{t HatHirMiuMol Thm 24 for metric spaces}, and hence the mentioned result implies that $U_{w} (\Delta(v)^*) = \Delta(v).$ Since $\|w -\Delta(v)\|<2,$ \cite[Lemma 2.3]{CuPe20b} guarantees that $w=\Delta(v).$ Having in mind that $U_w ( \Delta(u)^* ) = \Delta\left( U_v (u^*) \right)$ we get $\Delta(U_{v} (u^{*}))   = U_{\Delta(v)} (\Delta(u)^*),$ as desired.
\end{proof}

The next technical result is a refinement of \cite[Lemma 3.3]{CuPe20b}, which is in fact a Jordan version of \cite[Lemma 7]{HatMol2012}. The main difference is that in the present result we have a surjective isometry between the principal components of the sets of unitaries in two unital JB$^*$-algebras, while in \cite[Lemma 3.3]{CuPe20b} we find a surjective isometry between the whole sets of unitaries of two unital JB$^*$-algebras. Since the elements in the principal component is a self-adjoint quadratic subset of the set of unitaries (cf. Proposition \ref{p self adjoint quadratic subset princ comp}), the arguments in the proof of \cite[Lemma 3.3]{CuPe20b} are also valid in our setting. So, the details will be omitted.

\begin{lemma}\label{l HM Lemma 7 Jordan} Let $M$ and $N$ be two unital JB$^*$-algebras. Let $\{u_{k}: 0\leq k\leq 2^n\}$ be a subset of $\mathcal{U}^0(M)$ {\rm(}with $n \in \mathbb{N}${\rm)} and let $\Phi: \mathcal{U}^0(M)\to \mathcal{U}^0(N)$ be a mapping such that $$U_{u_{k+1}} (u_k^*) = u_{k+2}, \hbox{ and } \Phi (U_{u_{k+1}} (u_k^*)) = U_{\Phi(u_{k+1})} (\Phi (u_k)^*),$$ for all $0\leq k\leq 2^n -2$. Then $U_{u_{2^{n-1}}} (u_0^*) = u_{2^n}$ and $$\Phi\left(U_{u_{2^{n-1}}} (u_0^*) \right) = U_{\Phi(u_{2^{n-1}})} \Phi(u_0)^*. $$
\end{lemma}

Before addressing our main goal, we shall explore the behavior of those unital surjective isometries over $U$-products of exponentials of self-adjoint elements. The next result is the technical core of the arguments.

\begin{proposition}\label{p unital SI over U-products of exp} For each natural $n$ the identity $$ \Delta_0(U_{e^{i h_n}}\cdots U_{e^{i h_1}}(e^{i h_0}))=U_{\Delta_0(e^{i h_n})}\cdots U_{\Delta_0(e^{i h_1})}(\Delta(e^{i h_0})),$$ holds for every couple of unital JB$^*$-algebras $M$ and $N$, every  unital {\rm(}i.e., $\Delta_0(\11b{M})=\11b{N}${\rm)} surjective isometry $\Delta_0: \mathcal{U}^0 (M)\to \mathcal{U}^0 (N)$ and every $h_n,\dots,h_1, h_0\in M_{sa}$. Furthermore, if we fix a couple of unital JB$^*$-algebras $M$ and $N$ and we denote by Iso$_{1}(\mathcal{U}^0 (M), \mathcal{U}^0 (N))$ the set of all surjective unital isometries between $\mathcal{U}^0 (M)$ and $\mathcal{U}^0 (N)$, then there exists a mapping $k:M_{sa}\times \hbox{Iso}_{1} (\mathcal{U}^0 (M), \mathcal{U}^0 (N))\to N_{sa}$ satisfying the following properties:\begin{enumerate}[$(a)$]\item If we fix $\Delta_0\in \hbox{Iso}_{1} (\mathcal{U}^0 (M), \mathcal{U}^0 (N)),$ the mapping $k(\cdot, \Delta_0): M_{sa}\to N_{sa}$ is a surjective linear isometry;
\item $k(\cdot, \Delta_0)^{-1} = k(\cdot, \Delta_0^{-1}),$ for every $\Delta_0\in \hbox{Iso}_{1} (\mathcal{U}^0 (M), \mathcal{U}^0 (N));$
\item $\Delta_0(e^{i t h } )= e^{i t k(h, \Delta_0)}$ for all $\Delta_0\in \hbox{Iso}_{1} (\mathcal{U}^0 (M), \mathcal{U}^0 (N))$, $h\in M_{sa}$ and $t\in \mathbb{R}$.
\end{enumerate}

\end{proposition}

\begin{proof} Concerning the first identity, we observe that it suffices to prove that for each natural $n$ the identity \begin{equation}\label{eq ideneity at 1} \Delta_0(U_{e^{i h_n}}\cdots U_{e^{i h_1}}(\11b{M}))=U_{\Delta_0(e^{i h_n})}\cdots U_{\Delta_0(e^{i h_1})}(\11b{N})
 \end{equation} holds for every $\Delta_0$ and $h_n,\dots,h_1\in M_{sa}$ as in the hypotheses. Namely, if $h_0\in M_{sa}$ and the previous identity \eqref{eq ideneity at 1} is true we have $$\begin{aligned}\Delta_0(U_{e^{i h_n}}\cdots U_{e^{i h_1}}(e^{i h_0}))&= \Delta_0(U_{e^{i h_n}}\cdots U_{e^{i h_1}} U_{e^{i \frac{h_0}{2}}}(\11b{M}) ) \\
=U_{\Delta_0(e^{i h_n})}\cdots U_{\Delta_0(e^{i h_1})} U_{\Delta_0(e^{i \frac{h_0}{2}})} (\11b{N}) &= U_{\Delta_0(e^{i h_n})}\cdots U_{\Delta_0(e^{i h_1})} \Delta_0 (U_{e^{i \frac{h_0}{2}}} (\11b{M}) ) \\
&=U_{\Delta_0(e^{i h_n})}\cdots U_{\Delta_0(e^{i h_1})}(\Delta_0(e^{i h_0})),
\end{aligned} $$ where at the penultimate equality we applied \eqref{eq ideneity at 1}.\smallskip

We shall prove \eqref{eq ideneity at 1} by induction on $n$. Let us first assume that $n=1$, and fix $\Delta_0$ and $h_1\in M_{sa}$ as in the hypotheses. Fix an arbitrary $h_0\in M_{sa}$ and consider the continuous mapping $E_{h_0}: \mathbb{R}\to \mathcal{U}^0 (M)$, $E_{h_0}(t)=e^{i t h_0}$.  We begin our argument by proving that the identity
\begin{equation}\label{eq n1 claim1}
\Delta_0( U_{v_t}(v_s))=
U_{\Delta_0(v_t)}(\Delta_0( v_s^* )^*)
\end{equation} holds for every $v_s=e^{i s h_0}=E_{h_0} (s)$ and $v_t=e^{i t h_0}=E_{h_0}(t)$ in $E_{h_0}(\RR)$.\smallskip

Namely, choose a positive integer $m$ such that \begin{equation}\label{eq m} e^{\frac{\| i(s+t) h_0 \|_{_{M}}}{2^m}}-1 < \frac{1}{2},\hbox{ and hence } \left\| e^{\frac{i(s+t) h_0}{2^m}}-\11b{M} \right\|_{_{M}}\leq e^{\frac{\| i(s+t) h_0 \|_{_{M}}}{2^m}}-1 < \frac{1}{2}.
\end{equation}
	
Consider now the family $\{ u_l : 0\leq l \leq 2^{m+1}\}\subseteq \mathcal{U}^0(M)$, where $$u_l = v_s^* \circ e^{\frac{i l (s+t) h_0}{2^m}}= e^{-i s h_0}\circ e^{\frac{i l (s+t) h_0}{2^m}} =\hbox{\eqref{eq exp commute}}=e^{\left(-i s +\frac{i l (s+t)}{2^m}\right)  h_0}, \quad 0\leq l \leq 2^{m+1}.$$
	
Having in mind \eqref{eq exp commute} it is not hard to see that
	\begin{equation}\label{eq Uuk+1 is uk+2} \begin{aligned} u_0& =v_s^*,\  u_{2^m}=v_t,\  u_{2^{m+1}}=U_{v_t} (v_s), \\
	 U_{u_{l+1}} ({u_l}^*) & =U_{u_{l+1}}(u_l^{-1})= u_{l+2},\hbox{ for any $0\leq l \leq 2^{m+1}-2$,}\\
\hbox{ and }  &\|u_{l+1} - u_l \| = \left\| {v_s}^*\circ e^{\frac{i(l+1)(s+t) h_0}{2^m}}- {v_s}^*\circ e^{\frac{i l (s+t) h_0}{2^m}} \right\| \\
	&\leq  \|{v_s}^*\| \ \left\| e^{\frac{i l (s+t) h_0}{2^m}}\circ e^{\frac{i(s+t) h_0}{2^m}} - e^{\frac{i l (s+t) h_0}{2^m}} \right\|\\
	& \leq \left\|e^{\frac{i l (s+t) h_0}{2^m}}\right\|\  \left\| e^{\frac{i (s+t) h_0}{2^m}}-\11b{M} \right\| \\
&\leq e^{\frac{\| i(s+t) h_0 \|}{2^m}}-1 < \frac{1}{2},\hbox{ for all } 0\leq l \leq 2^{m+1}-1.
\end{aligned}
	\end{equation}

An application of Theorem \ref{t Delta preserves inverted triple products in U0M}$(3)$ to $\Delta_0$ and the unitaries $u_{l+1}$ and $u_l$ gives
	\begin{equation}\label{eq application of Th2.9 3}
	\Delta_0(U_{u_{l+1}} ({u_l}^*)) = U_{\Delta_0(u_{l+1})} (\Delta_0(u_{l})^*), \hbox{  for every $l\in \{0,\dots, 2^{m+1}-1\}$.}
	\end{equation}

It follows from \eqref{eq Uuk+1 is uk+2} and \eqref{eq application of Th2.9 3} that we are in a position to apply Lemma \ref{l HM Lemma 7 Jordan} with $\Phi = \Delta_0$ and $n= m+1$ to deduce that
$$\begin{aligned}
	\Delta_0 (U_{v_t} ({v_s}) )&= \Delta_0 (U_{v_t} (({v_s}^*)^*)) =\Delta_0(U_{u_{2^m}} ({u_0}^*)) \\
&= U_{\Delta_0(u_{2^m})} (\Delta_0(u_{0})^*) =U_{\Delta_0 (v_t)} (\Delta_0({v_s}^*)^*),
\end{aligned}$$  which concludes the proof of \eqref{eq n1 claim1}.\smallskip

Let us see some consequences derived from \eqref{eq n1 claim1}. If in this identity we take $t=1$, $s=0$ and $h_0=h_1\in M_{sa}$ we get
\begin{equation}\label{eq Deltaa1 to the square}
\Delta_0({e^{i 2 h_1}}) =\Delta_0( U_{e^{i h_1}} (\11b{M}))= U_{\Delta_0(e^{i h_1})}(\11b{N})= \Delta_0(e^{i h_1})^2,
 \end{equation} which, in particular, proves \eqref{eq ideneity at 1} for $n=1$.\smallskip

Second, by replacing $t=0$ and $h_0=h_1\in M_{sa}$ in \eqref{eq n1 claim1} we get
	\begin{equation}\label{eq Delta0 adjoint} \Delta_0(e^{i s h_1}) =\Delta_0(e^{-i s h_1})^*, \ \ \forall h_1\in M_{sa}, \ \forall s\in \RR.\end{equation}

The most important consequence derived from \eqref{eq n1 claim1} (and \eqref{eq Delta0 adjoint}) assures that, for each $h_0\in M_{sa}$, the mapping $$\mathbb{R}\to \mathcal{U}^0 (N), \ t\mapsto \Delta_0(e^{i t h_0})$$ is a uniformly continuous one-parameter family of unitaries satisfying the hypotheses of the Jordan version of Stone's one-parameter theorem in Theorem \ref{t Jordan unitary groups version of Stone's theorem} \cite[Theorem 3.1]{CuPe20b}. Namely, $$\begin{aligned}U_{ \Delta_0(e^{i t h_0}) } ( \Delta_0(e^{ i s h_0}) ) &= U_{ \Delta_0(e^{i t h_0}) } ( \Delta_0(e^{- i s h_0})^* ) =  \Delta_0( U_{e^{i t h_0}}(e^{i s h_0}))\\
&=  \Delta_0( e^{i(2 t + s) h_0}).
 \end{aligned}$$ Therefore there exists a unique $k(h_0, \Delta_0)\in N_{sa},$ depending on $h_0$ and $\Delta_0,$ satisfying \begin{equation}\label{eq application Stone one parameter thm}  \Delta_0(e^{i t h_0}) = e^{i t  k(h_0, \Delta_0)}, \hbox{ for all } t\in \mathbb{R}.
\end{equation}

We observe that the mapping $k(\cdot, \Delta_0): M_{sa} \to N_{sa}$ is well defined, and moreover, by \eqref{eq application Stone one parameter thm}, $\11b{N}= \Delta_0(e^{i t 0}) = e^{i t  k(0, \Delta_0)},$ for all $t\in \mathbb{R},$ which assures that $k(0, \Delta_0)=0$. Let us prove that this mapping is a surjective isometry, and as a consequence of the Mazur--Ulam theorem, it must be a surjective linear isometry.\smallskip

Similarly, by considering the surjective linear isometry $\Delta_0^{-1}: \mathcal{U}^0 (N) \to \mathcal{U}^{0} (M)$, the above arguments prove that  for each $k_0\in N_{sa}$ there exists a unique $h(k_0, \Delta_0^{-1})\in M_{sa},$ depending on $k_0$ and $\Delta_0^{-1},$ satisfying \begin{equation}\label{eq application Stone one parameter thm Delta inverse}  \Delta_0^{-1} (e^{i t k_0}) = e^{i t  h(k_0, \Delta_0^{-1})}, \hbox{ for all } t\in \mathbb{R}.
\end{equation} As before, the mapping $h(\cdot, \Delta_0^{-1}): N_{sa} \to M_{sa}$ is well defined. By combining \eqref{eq application Stone one parameter thm} and \eqref{eq application Stone one parameter thm Delta inverse} we get $$\begin{aligned}e^{i t h_0} &= \Delta_0^{-1} \Delta_0(e^{i t h_0}) = \Delta_0^{-1} \left( e^{i t  k(h_0, \Delta_0)} \right) = e^{i t  h(k(h_0, \Delta_0), \Delta_0^{-1})},\ \hbox{ and }\\
e^{i t k_0} &= \Delta_0 \Delta_0^{-1} (e^{i t k_0}) = \Delta_0 \left( e^{i t  h(k_0, \Delta_0^{-1})} \right) = e^{i t  k(h(k_0, \Delta_0^{-1}), \Delta_0)},
\end{aligned}  $$ for all $t\in \mathbb{R}$, $h_0\in M_{sa}$, $k_0\in N_{sa}$. Taking a simple derivative at $t=0$ we arrive to the identities $$h(k(h_0, \Delta_0), \Delta_0^{-1}) = h_0, \hbox{ and } k(h(k_0, \Delta_0^{-1}), \Delta_0) = k_0, $$ for all $h_0\in M_{sa}$, $k_0\in N_{sa}$, that is, $k(\cdot, \Delta_0)$ and $h(\cdot, \Delta_0^{-1})$ are bijections with $h(\cdot, \Delta_0^{-1}) = k(\cdot, \Delta_0)^{-1}$.\smallskip

We shall next prove that $k(\cdot, \Delta_0): M_{sa} \to N_{sa}$ is an isometry. This is a standard procedure already in \cite{HatMol2014,Hatori14,CuPe20b}. Given $h_0,h_0'\in M_{sa}$ and a real number $t$ it is well known that
$$\|\cdot\|\hbox{-}\lim_{t\to 0^+} \frac{e^{i t h_0}-e^{ith_0'}}{t}= \|\cdot\|\hbox{-}\lim_{t\to 0^+} \frac{e^{i t h_0}-\11b{M}}{t}-\frac{e^{i t h_0'}-\11b{M}}{t}  = ih_0-ih_0'.$$
It follows from \eqref{eq application Stone one parameter thm} that $$ \|e^{it k(h_0, \Delta_0)} - e^{it k(h_0', \Delta_0)} \|_{_{N}}=\|\Delta_0 (e^{ith_0})-\Delta_0 (e^{ith_0'}) \|_{_{N}}=\| e^{i t h_0}-e^{it h_0'} \|_{_{M}}. $$ Therefore, we deduce from the uniqueness of the limits above that  $$\|k(h_0, \Delta_0) - k(h_0', \Delta_0) \|_{_{N}}=\|h_0-h_0' \|_{_{M}},$$ witnessing that $k(\cdot, \Delta_0) :M_{sa}\to N_{sa}$ is a surjective isometry, and hence a surjective linear isometry by the Mazur--Ulam theorem. We have shown the following, for each unital surjective isometry $\Delta_0: \mathcal{U}^0 (M) \to \mathcal{U}^0 (N)$ \begin{equation}\label{eq the log associated with the exponential of Delta0} \hbox{the mapping } k(\cdot, \Delta_0): M_{sa} \to N_{sa} \hbox{ is a surjective linear isometry.}
\end{equation}

We shall next prove the following identity
\begin{equation}\label{eq n1 with 2 exponentials}
\Delta_0( U_{e^{i h_1}}(e^{i h_0}))= U_{\Delta_0 (e^{i h_1})} (\Delta_0(e^{i h_0})), \hbox{ for all } h_1,h_0\in M_{sa}.
\end{equation}

Namely, fix $h_1\in M_{sa}$ and consider the surjective isometry $\Delta_0^{h_1} :=U_{\Delta_0 (e^{i h_1})^*} \Delta_0 U_{e^{i h_1}} : \mathcal{U}^{0} (M) \to \mathcal{U}^{0} (N),$ which is well defined by Theorem \ref{t principal component as product of exp}. By \eqref{eq Deltaa1 to the square} we obtain
$$\begin{aligned} \Delta_0^{h_1} (\11b{M}) &= U_{\Delta_0(e^{i h_1})^*} \Delta_0 U_{e^{i h_1}} (\11b{M}) = U_{\Delta_0(e^{i h_1})^*} \Delta_0 (e^{ 2 i h_1}) =  U_{\Delta_0(e^{i h_1})^*} (\Delta_0 (e^{ i h_1})^2) \\
&= U_{\Delta_0(e^{i h_1})^*} U_{\Delta_0 (e^{ i h_1})} (\11b{M})=  \11b{N},
\end{aligned} $$ witnessing that the mapping $\Delta_0^{h_1}  = U_{\Delta_0(e^{i h_1})^*} \Delta_0 U_{e^{i h_1}}$ is unital.\smallskip

Let $k(\cdot, \Delta_0), k(\cdot, \Delta_0^{h_1}): M_{sa} \to N_{sa}$ be the surjective linear isometries given by \eqref{eq application Stone one parameter thm} and \eqref{eq the log associated with the exponential of Delta0} for $\Delta_0$ and $\Delta_0^{h_1}$, respectively.\smallskip

Choose, by continuity, $\varepsilon>0$ such that
$$ \|e^{i h_1}-e^{-i t h}\|<\frac12, \ \forall t\in \RR, \ h\in M_{sa} \hbox{ with } |t-1|<\varepsilon, \hbox{ and } \|h + h_1\|<\varepsilon.$$
By applying Theorem \ref{t Delta preserves inverted triple products in U0M}$(3)$ to $\Delta_0$ and the unitaries $e^{i h_1}, e^{-it h}\in \mathcal{U}^0(M)$ with $t\in \RR$ and $h\in M_{sa}$ satisfying $\|h+h_1\|<\varepsilon$ and $|t-1|<\varepsilon$ we conclude that
	$$\Delta_0(U_{e^{i h_1}} (e^{it h})) =U_{\Delta_0( e^{i h_1} )} ((\Delta_0(e^{-it h}))^*)=\eqref{eq Delta0 adjoint}=U_{\Delta_0(e^{i h_1})} (\Delta_0(e^{it h})).$$
The previous identity implies, by the definition of $k(\cdot, \Delta_0)$ and $k(\cdot, \Delta_0^{h_1})$, that $$e^{i t  k(h , \Delta_0)} = e^{i t  k(h, \Delta_0^{h_1})},\  \forall t\in \RR, \ h\in M_{sa} \hbox{ with } |t-1|<\varepsilon, \hbox{ and } \|h + h_1\|<\varepsilon.$$ If we fix an arbitrary $h\in M_{sa}$ with $\|h + h_1\|<\varepsilon$ and we take a simple derivative at $t=1$ in the above equality we obtain $$ k(h , \Delta_0)= k(h, \Delta_0^{h_1}), \hbox{ for all } h\in M_{sa} \hbox{ with } \|h + h_1\|<\varepsilon,$$ and the linearity of $k(\cdot, \Delta_0),$ $k(\cdot, \Delta_0^{h_1})$ implies that $k(\cdot, \Delta_0)=k(\cdot, \Delta_0^{h_1})$.\smallskip

Now, it follows from \eqref{eq application Stone one parameter thm} that $$ \Delta_0(e^{i t h_0}) = e^{i t  k(h_0, \Delta_0)}   = e^{i t  k(h_0, \Delta_0^{h_1})} = \Delta_0^{h_1} (e^{i t h_0}) = U_{\Delta_0 (e^{i h_1})^*} \Delta_0 U_{e^{i h_1}} (e^{i t h_0}) ,$$ for all real $t$ and all $h_0\in M_{sa}$, which concludes the proof of \eqref{eq n1 with 2 exponentials}.\smallskip

Let us now proceed with the final part in the induction argument to get  \eqref{eq ideneity at 1}. Suppose by the induction hypothesis that \eqref{eq ideneity at 1} holds for every unital surjective isometry $\Delta_0$ and $h_n,\dots,h_1\in M_{sa}$ as in the hypotheses. Take an arbitrary unital surjective isometry $\Delta_1 : \mathcal{U}^0 (M)\to \mathcal{U}^0 (N)$ and $h_{n+1},h_n,\dots,h_1\in M_{sa}$. Set $w_{0} = e^{2 i h_{n+1}}$ and $\widetilde{w}_{0} = \Delta_1(e^{2 i h_{n+1}})\in \mathcal{U}^0 (N).$ We consider the unital JB$^*$-algebras given by the isotopes $M(w_{0})$ and $N(\widetilde{w}_{0})$ with the corresponding Jordan products $\circ_{w_{0}}$ and $\circ_{\widetilde{w}_{0}}$ and involutions $*_{w_{0}}$ and $*_{\widetilde{w}_{0}}$, respectively. Since $w_{0}\in \mathcal{U}^0 (M)$ and $\widetilde{w}_{0} \in \mathcal{U}^0 (N),$ it follows from \eqref{eq princ components of isotopes coincide} that $$\mathcal{U}^0(M(w_{0}))=\mathcal{U}^0(M), \hbox{ and } \mathcal{U}^0(N(\widetilde{w}_{0}))=\mathcal{U}^0(N),$$ and thus $\widetilde{\Delta}_1 : \mathcal{U}^0(M(w_{0}))=\mathcal{U}^0(M)\to \mathcal{U}^0(N(\widetilde{w}_{0}))=\mathcal{U}^0(N),$ $\widetilde{\Delta}_1 (a) = {\Delta}_1(a)$ is a unital (i.e., $\Delta_1(w_{0}) = \widetilde{w}_{0}$) surjective isometry.\smallskip

The mapping $U_{e^{i h_{n+1}}} : M\to M(w_{0})$ is a unital surjective linear isometry, and hence a Jordan $^*$-isomorphism between these two unital JB$^*$-algebras (cf. Remark \ref{r isotopes remark} or \cite{WriYou78}). Therefore $\hat{h}_k:= U_{e^{i h_{n+1}}} (h_k) \in (M(w_{0}))_{sa}$ for all $1\leq k \leq n$.\smallskip

By applying the induction hypothesis to $\widetilde{\Delta}_1 : \mathcal{U}^0(M(w_{0}))\to \mathcal{U}^0(N(\widetilde{w}_{0})),$ and $\hat{h}_1, \ldots, \hat{h}_n$ we obtain
\begin{equation}\label{eq induction hypothesis applied to tildeDelta1}\begin{aligned} &{\Delta}_1 U^{(w_{0})}_{\exp_{w_{0}} (i \hat{h}_{n})} \ldots U^{(w_{0})}_{\exp_{w_{0}} (i \hat{h}_{1})} (w_{0}) = \widetilde{\Delta}_1 U^{(w_{0})}_{\exp_{w_{0}} (i \hat{h}_{n})} \ldots U^{(w_{0})}_{\exp_{w_{0}} (i \hat{h}_{1})} (w_{0}) \\
&= U^{(\widetilde{w}_{0})}_{ \widetilde{\Delta}_1(\exp_{w_{0}} (i \hat{h}_{n}))} \ldots U^{(\widetilde{w}_{0})}_{ \widetilde{\Delta}_1(\exp_{w_{0}} (i \hat{h}_{1}))} (\widetilde{w}_{0}) .
\end{aligned}
\end{equation}

Let us carefully decode the identity in \eqref{eq induction hypothesis applied to tildeDelta1}. First by  \eqref{eq quadratic operation on isotopes} and the fact that $M(w_0) = M_{(w_0^*)},$ for any $v_1,v_2\in \mathcal{U}^{0} (M)$ we have \begin{equation}\label{eq 030521 a}
U^{(w_0)}_{v_1} (v_2) = U_{v_1} U_{w_0^*} (v_2)  = U_{v_1} U_{e^{-2 i h_1}} (v_2).
  \end{equation}

Furthermore, since, as before, $\exp_{w_0}$ denotes the exponential in the unital JB$^*$-algebra $M(w_0)$, it follows from the fact that $U_{e^{i h_{n+1}}} : M\to M(w_0)$ is a Jordan $^*$-isomorphism that \begin{equation}\label{eq exponentials of hathk} \exp_{w_0} (i \hat{h}_k)= \exp_{w_0} (i U_{e^{i h_{n+1}}} (h_k)) = U_{e^{i h_{n+1}}} (e^{ih_k}), \hbox{ for all } 1\leq k\leq n,
 \end{equation} and by the fundamental identity \eqref{eq fundamental identity UaUbUa} we deduce
\begin{equation}\label{eq exponentials and U} \begin{aligned}
 U_{\exp_{w_0} (i \hat{h}_k)}^{(w_0)} & = U_{\exp_{w_0} (i \hat{h}_k)} U_{{w_0}^*}  = U_{U_{e^{i h_{n+1}}} (e^{ih_k})} U_{{w_0}^*} \\
 &= U_{e^{i h_{n+1}}} U_{e^{ih_k}} U_{e^{i h_{n+1}}} U_{{e^{- 2 i h_{n+1}}}} = U_{e^{i h_{n+1}}} U_{e^{ih_k}} U_{e^{-i h_{n+1}}},
\end{aligned}
\end{equation} where in the last step we applied that $e^{i h_{n+1}}$ and $e^{-2i h_{n+1}}$ operator commute. \smallskip

Having in mind \eqref{eq 030521 a} and \eqref{eq exponentials and U}, the left hand side of \eqref{eq induction hypothesis applied to tildeDelta1} writes in the form \begin{equation}\label{eq left hand side of 25} \begin{aligned} & {\Delta}_1 U^{(w_{0})}_{\exp_{w_{0}} (i \hat{h}_{n})} \ldots U^{(w_{0})}_{\exp_{w_{0}} (i \hat{h}_{1})} (w_{0}) \\
&=  {\Delta}_1 U_{e^{i h_{n+1}}} U_{e^{ih_n}} U_{e^{-i h_{n+1}}} \ldots U_{e^{i h_{n+1}}} U_{e^{ih_1}} U_{e^{-i h_{n+1}}}  (e^{2i h_{n+1}}) \\
&= {\Delta}_1 U_{e^{i h_{n+1}}} U_{e^{ih_n}}  U_{e^{i h_{n-1}}}  \ldots U_{e^{i h_{2}}} U_{e^{ih_1}} (\11b{M}).
\end{aligned}
\end{equation}

To deal with the right hand side of \eqref{eq induction hypothesis applied to tildeDelta1} we observe that, by \eqref{eq exponentials of hathk}, for $1\leq k\leq n$ we have
$$\begin{aligned}
U^{(\widetilde{w}_{0})}_{ \widetilde{\Delta}_1(\exp_{w_{0}} (i \hat{h}_{k}))} &=  U^{(\widetilde{w}_{0})}_{{\Delta}_1(\exp_{w_{0}} (i \hat{h}_{k}))} = U_{{\Delta}_1(U_{e^{i h_{n+1}}} (e^{ih_k}))} U_{\widetilde{w}_{0}^*} \\
&= \hbox{(by \eqref{eq n1 with 2 exponentials})} = U_{U_{{\Delta}_1(e^{i h_{n+1}})} {\Delta}_1 (e^{ih_k})} U_{\widetilde{w}_{0}^*} = \hbox{(by \eqref{eq fundamental identity UaUbUa})} \\
&= U_{{\Delta}_1(e^{i h_{n+1}})} U_{{\Delta}_1 (e^{ih_k})} U_{{\Delta}_1(e^{i h_{n+1}})} U_{\Delta_1(e^{2 i h_{n+1}})^*}\\
&= \hbox{(by \eqref{eq Deltaa1 to the square})} = U_{{\Delta}_1(e^{i h_{n+1}})} U_{{\Delta}_1 (e^{ih_k})} U_{{\Delta}_1(e^{i h_{n+1}})} U_{(\Delta_1(e^{ i h_{n+1}})^2)^*}\\
& = U_{{\Delta}_1(e^{i h_{n+1}})} U_{{\Delta}_1 (e^{ih_k})} U_{{\Delta}_1(e^{i h_{n+1}})^*} \\
&= U_{{\Delta}_1(e^{i h_{n+1}})} U_{{\Delta}_1 (e^{ih_k})} U_{{\Delta}_1(e^{i h_{n+1}})}^{-1},
\end{aligned} $$ and hence the right hand side of \eqref{eq induction hypothesis applied to tildeDelta1} writes in the form
\begin{equation}\label{eq right han side simplified}\begin{aligned}
& U^{(\widetilde{w}_{0})}_{ \widetilde{\Delta}_1(\exp_{w_{0}} (i \hat{h}_{n}))} \ldots U^{(\widetilde{w}_{0})}_{ \widetilde{\Delta}_1(\exp_{w_{0}} (i \hat{h}_{1}))} (\widetilde{w}_{0}) \\
&= U_{{\Delta}_1(e^{i h_{n+1}})} U_{{\Delta}_1 (e^{ih_n})} \ldots  U_{{\Delta}_1 (e^{ih_1})} U_{{\Delta}_1(e^{i h_{n+1}})}^{-1} (\Delta_1(e^{2 i h_{n+1}}))\\
&= \hbox{(by \eqref{eq Deltaa1 to the square})} = U_{{\Delta}_1(e^{i h_{n+1}})} U_{{\Delta}_1 (e^{ih_n})} \ldots  U_{{\Delta}_1 (e^{ih_1})} U_{{\Delta}_1(e^{i h_{n+1}})}^{-1} (\Delta_1(e^{ i h_{n+1}})^2)\\
& =  U_{{\Delta}_1(e^{i h_{n+1}})} U_{{\Delta}_1 (e^{ih_n})} \ldots  U_{{\Delta}_1 (e^{ih_1})} (\11b{N}). \end{aligned}
\end{equation}

Finally, \eqref{eq induction hypothesis applied to tildeDelta1}, \eqref{eq left hand side of 25} and \eqref{eq right han side simplified} give $$ {\Delta}_1 U_{e^{i h_{n+1}}} U_{e^{ih_n}}   \ldots U_{e^{i h_{2}}} U_{e^{ih_1}} (\11b{M})=  U_{{\Delta}_1(e^{i h_{n+1}})} U_{{\Delta}_1 (e^{ih_n})} \ldots  U_{{\Delta}_1 (e^{ih_1})} (\11b{N}).$$ By induction, we have \eqref{eq ideneity at 1} for every $\Delta_0,$ $n\in \mathbb{N}$, and $h_n, \ldots, h_1 \in M_{sa}$.
\end{proof}	

Let $p$ and $q$ be two projections in a JB$^*$-algebra. We briefly recall that $p$ and $q$ are called \emph{orthogonal} if $p \circ q = 0$. Actually the notion of orthogonality has been studied and developed for general elements in a JB$^*$-triple (cf. \cite{BurFerGarMarPe}). Elements $a,b$ in a JB$^*$-triple $E$ are called \emph{orthogonal} (written $a\perp b$) if $L(a,b) =0$. It is known that $a\perp b$ $\Leftrightarrow$ $\J aab =0$ $\Leftrightarrow$ $\{b,b,a\}=0$ $\Leftrightarrow$ $b\perp a$ (see \cite[Lemma 1]{BurFerGarMarPe}).

\begin{theorem}\label{t surj isom principal components}
Let $\Delta: \mathcal{U}^0 (M)\to \mathcal{U}^0 (N)$ be a surjective isometry between the principal components of two unital JB$^*$-algebras. Then there exist $k_1,\ldots,k_n\in N_{sa}$, a central projection $p\in N$ and a Jordan $^*$-isomorphism $\Phi:M\to N$ such that $$\begin{aligned}\Delta(u) &= p\circ U_{e^{i k_n}} \ldots  U_{e^{i k_1}}  \Phi (u) + (\11b{N}-p)\circ \left( U_{e^{-i k_n}} \ldots  U_{e^{-i k_1}} \Phi(u)\right)^*,
\end{aligned}$$ for all $u\in \mathcal{U}^0(M)$. Consequently, $M$ and $N$ are Jordan $^*$-isomorphic, and there exists a surjective real linear isometry {\rm(}i.e., a real linear triple isomorphism{\rm)} from $M$ onto $N$ whose restriction to $\mathcal{U}^0 (M)$ is $\Delta$.
\end{theorem}

\begin{proof} Since $\Delta(\11b{M})\in \mathcal{U}^0 (N)$, Theorem \ref{t principal component as product of exp} implies the existence of $k_1,\ldots,k_n\in N_{sa}$ such that $ \Delta(\11b{M}) = U_{e^{i k_n}} \ldots  U_{e^{i k_1}} (\11b{N})$. It follows from Proposition \ref{p self adjoint quadratic subset princ comp} that the mapping $\Delta_0 = U_{e^{-i k_1}} \ldots  U_{e^{-i k_n}} \Delta : \mathcal{U}^0 (M)\to \mathcal{U}^0 (N)$ is a surjective isometry, and by construction $\Delta_0 (\11b{M}) = \11b{N}$.\smallskip
	
Let $k(\cdot,\Delta_0): M_{sa}\to N_{sa}$ be the surjective linear isometry given by Proposition \ref{p unital SI over U-products of exp}. \smallskip
	
The self-adjoint part of any JB$^*$-algebra is a JB-algebra. Thus, $k(\cdot,\Delta_0) : M_{sa}\to N_{sa}$ is a surjective linear isometry between JB-algebras. Theorem 1.4 and Corollary 1.11 in \cite{IsRo95} guarantee the existence of a central symmetry $k(\11b{M},\Delta_0)$ in $N_{sa}$ and a Jordan $^*$-isomorphism $\Phi:M\to N$ such that \begin{equation}\label{eq f description}
	k(h,\Delta_0) = k(\11b{M},\Delta_0) \circ \Phi(h),
	\end{equation} for every $h\in M_{sa}$. 
\smallskip
	
The rest of the proof follows similar arguments to those in the proof of \cite[Theorem 3.4]{CuPe20b}. That is, by construction there exists a central projection $p$ in $N$ such that $k(\11b{M},\Delta_0)= 2p -\11b{N}= p - (\11b{N}-p)$, where $p$ and $(\11b{N}-p)$ clearly are orthogonal projections in $N$, and for any $n>0$, $$\left( 2p-\11b{N}\right)^n= \left( p-(\11b{N}-p)\right)^n= p + (-1)^n(\11b{N}-p).$$

The mapping $\Delta_0$ can be expressed in terms of $p$ and $\Phi$ by some simple computations. Namely, given an arbitrary $h\in M_{sa}$, we have
$$\begin{aligned}
	\Delta_0(e^{i h})&=e^{i k(h,\Delta_0)}=e^{i k(\11b{M},\Delta_0)\circ \Phi(h) }=e^{i(2p-\11b{N})\circ \Phi(h) }= \sum_{n=0}^{\infty}\frac{\left( i  (2p-\11b{N})\circ \Phi(h)\right)^n }{n!}.
	\end{aligned}$$
We now make use of the properties of $\Phi$ as Jordan $^*$-isomorphism and the expression of $\left( 2p-\11b{N}\right)^n$ given above. It is worth noting that $k(\11b{M},\Delta_0)$ (and hence $p$) is central, and so it operator commutes with any element in $N$. Therefore,
	\begin{equation}\label{eq Delta0 in eix}\begin{aligned}
	\Delta_0(e^{ih})&=\sum_{n=0}^{\infty}\frac{ i^n  (2p-\11b{N})^n\circ \Phi(h)^n }{n!} =\sum_{n=0}^{\infty}\frac{ i^n  \left( p + (-1)^n(\11b{N}-p)\right) \circ \Phi(h^n) }{n!} \\
	&= \sum_{n=0}^{\infty}\frac{ i^n p \circ \Phi(h^n) }{n!} + \sum_{n=0}^{\infty}\frac{ i^n (-1)^n(\11b{N}-p) \circ \Phi(h^n) }{n!} \\
	&=p \circ \sum_{n=0}^{\infty}\frac{ i^n \Phi(h^n) }{n!} + (\11b{N}-p) \circ\sum_{n=0}^{\infty}\frac{ i^n (-1)^n \Phi(h^n) }{n!} \\
	&=p \circ \Phi\left( \sum_{n=0}^{\infty}\frac{ i^n h^n }{n!}\right)  + (\11b{N}-p) \circ\Phi\left( \sum_{n=0}^{\infty}\frac{ i^n (-1)^n h^n}{n!}\right) \\
	&= p \circ \Phi(e^{ih})  + (\11b{N}-p) \circ\Phi( e^{-ih})= p \circ \Phi(e^{ih})  + (\11b{N}-p) \circ\Phi( e^{ih})^*.
	\end{aligned}
\end{equation}
	
Let us fix an arbitrary element $u\in \mathcal{U}^0(M)$. By Theorem \ref{t principal component as product of exp}, there exist a natural $m$ and $h_1,\dots, h_m\in M_{sa}$ such that $u=U_{e^{i h_m}}\cdots U_{e^{i h_1}}(\11b{M})$. According to this expression, we can apply Proposition \ref{p unital SI over U-products of exp} and \eqref{eq Delta0 in eix} to compute the evaluation of $\Delta_0$ at the element $u$:
$$\begin{aligned}
\Delta_0(u)&=\Delta_0(U_{e^{i h_m}}\cdots U_{e^{i h_1}}(\11b{M})) = U_{\Delta_0(e^{i h_m})}\cdots U_{\Delta_0(e^{i h_1})}(\Delta_0(\11b{M})) \\
	&=U_{\left( p \circ \Phi(e^{i h_m})  + (\11b{N}-p) \circ\Phi( e^{i h_m})^* \right) }\cdots U_{\left( p \circ \Phi(e^{i h_1})  + (\11b{N}-p) \circ\Phi( e^{i h_1})^*\right) }(\11b{N})
	\end{aligned}$$

Since $p$ is a central projection in $N$, it is not hard to see that $$U_{p\circ a+(\11b{N}-p)\circ b} = U_{p\circ a} + U_{(\11b{N}-p)\circ b} =p\circ U_{ a} +(\11b{N}-p)\circ U_{ b},$$ for all $a,b\in N$. The orthogonality of $p$ and $\11b{N}-p$ can be now combined with the fact that $\Phi$ is a Jordan $^*$-isomorphism to conclude that
$$\begin{aligned}
\Delta_0(u)&= p\circ \Phi (U_{e^{i h_m}}\cdots U_{e^{i h_1}} (\11b{M})) + (\11b{N}-p)\circ \Phi(U_{e^{i h_m}}\cdots U_{e^{i h_1}}(\11b{M}))^*\\
	&= p\circ \Phi (u) + (\11b{N}-p)\circ \Phi(u)^*.
	\end{aligned}$$
	
Finally, by the definition of $\Delta_0$, we arrive at
	$$\begin{aligned}
		\Delta(u)&= U_{e^{i k_n}} \ldots  U_{e^{i k_1}} \Delta_0 (u)= U_{e^{i k_n}} \ldots  U_{e^{i k_1}} \left( p\circ \Phi (u) + (\11b{N}-p)\circ \Phi(u)^*\right),\\
&=   p\circ U_{e^{i k_n}} \ldots  U_{e^{i k_1}}  \Phi (u) + (\11b{N}-p)\circ \left( U_{e^{-i k_n}} \ldots  U_{e^{-i k_1}} \Phi(u)\right)^*,
	\end{aligned}$$ for all $u\in \mathcal{U}^0(M),$ which concludes the proof.
\end{proof}

We can now strengthen the conclusion in \cite[Corollary 3.8]{CuPe20b} where it is deduced that two unital JB$^*$-algebras $M$ and $N$ are isometrically isomorphic as (complex) Banach spaces if, and only if, they are isometrically isomorphic as real Banach spaces if, and only if, there exists a surjective isometry $\Delta: \mathcal{U}(M)\to \mathcal{U}(N).$

\begin{corollary}\label{c two unital JB*algebras are isomorphic iff their unitaries are isometric iff their principal components are isometric}
Let $M$ and $N$ be two unital JB$^*$-algebras. Then the following statements are equivalent:
\begin{enumerate}[$(a)$]
\item $M$ and $N$ are Jordan $^*$-isomorphic;
\item There exists a surjective isometry $\Delta: \mathcal{U}(M)\to \mathcal{U}(N)$ satisfying $\Delta (\11b{M})\in \mathcal{U}^0(N)$;
\item There exists a surjective isometry $\Delta: \mathcal{U}^0(M)\to \mathcal{U}^0(N).$
\end{enumerate}
\end{corollary}

\begin{proof} The implication $(a)\Leftrightarrow (b)$ is almost explicit in \cite[Corollary 3.8]{CuPe20b}, we include here an alternative argument.\smallskip

The implication $(a)\Rightarrow (b)$ is clear because every Jordan $^*$-isomorphism is an isometry and maps unitaries to unitaries.\smallskip

$(b)\Rightarrow (c)$ Suppose there exists a surjective isometry $\Delta: \mathcal{U}(M)\to \mathcal{U}(N)$ satisfying $\Delta (\11b{M})\in \mathcal{U}^0(N)$. Since every surjective isometry maps connected components to connected components, $\Delta(\mathcal{U}^0(M))$ is a connected component of $\mathcal{U}(N)$ and contains an element in $\mathcal{U}^0(N)$, the equality $\Delta(\mathcal{U}^0(M))= \mathcal{U}^0(N)$ holds, and thus $\Delta|_{\mathcal{U}^0(M)} : \mathcal{U}^0(M)\to \mathcal{U}^0(N)$ is a surjective isometry.\smallskip

Finally, the  implication $(c)\Rightarrow (a)$ is a consequence of Theorem \ref{t surj isom principal components}.
\end{proof}

If in Corollary \ref{c two unital JB*algebras are isomorphic iff their unitaries are isometric iff their principal components are isometric} the set $\mathcal{U} (M)$ is connected all the statements are also equivalent to the existence of a surjective isometry $\Delta: \mathcal{U}(M)\to \mathcal{U}(N)$.\smallskip

Let $M$ be a unital JB$^*$-algebra, let $\mathcal{U}^c (M)$ be a connected component of $\mathcal{U} (M)$, and let $w$ be an element in $\mathcal{U}^c (M)$. We shall be mainly interested in the case in which $\mathcal{U}^c (M) \neq \mathcal{U}^0 (M)$. We have already commented that $\mathcal{U}^c (M)$ is the principal component of the $w$-isotope $M(w)$ (cf. Remark \ref{r connected components which are not the pricipal one}). When particularized to the corresponding isotopes, the previous Theorem \ref{t surj isom principal components} can be re-stated as follows.

\begin{theorem}\label{t surj isom secondary components} Let $M$ and $N$ be unital JB$^*$-algebras. Let $\Delta: \mathcal{U}^c (M)\to \mathcal{U}^c (N)$ be a surjective isometry between two connected components of $M$ and $N$. Suppose $w_1\in \mathcal{U}^c (M)$ and $w_2\in \mathcal{U}^c (N)$. Then there exist $k_1,\ldots,k_n\in (N(w_2))_{sa}$, a central projection $p\in N(w_2)$ and a Jordan $^*$-isomorphism $\Phi:M(w_1)\to N(w_2)$ such that $$\begin{aligned}\Delta(u) &= p\circ_{w_2} U^{(w_2)}_{\exp_{w_2}({i k_n})} \ldots  U^{(w_2)}_{\exp_{w_2}({i k_1})}  \Phi (u) \\
&+ (w_2-p)\circ_{w_2} \left( U^{(w_2)}_{\exp_{w_2}({-i k_n})} \ldots  U^{(w_2)}_{\exp_{w_2}({-i k_1})} \Phi(u)\right)^{*_{w_2}},
\end{aligned}$$ for all $u\in \mathcal{U}^c (M)$. In particular, the unital JB$^*$-algebras $M(w_1)$ and $N(w_2)$ are Jordan $^*$-isomorphic, and there exists a surjective real linear isometry {\rm(}i.e., a real linear triple isomorphism{\rm)} from $M$ onto $N$ whose restriction to $\mathcal{U}^c (M)$ is $\Delta$.
\end{theorem}

We are now in position to exhibit additional, and deep, differences between the natural properties of the connected components of the unitary group in a unital C$^*$-algebra with respect to those in the unitary set of a unital JB$^*$-algebra.

\begin{remark}\label{r existence of unitaries with non isomorphic nor isometric connected components in the Jordan setting} R. Braun, W. Kaup and H. Upmeier show in \cite[Example 5.7]{BraKaUp78} that there exists a unital JB$^*$-algebra $M$ with unit $\11b{M}$ containing a unitary element $w$ for which we cannot find a surjective linear isometry $T: M\to M$ mapping $\11b{M}$ to $w$, that is, the group of all surjective linear isometries on $M$ --equivalently, triple automorphisms-- on $M$ is not transitive on $\mathcal{U} (M)$. It is further shown in \cite[Antitheorem 3.4.34]{Cabrera-Rodriguez-vol1} that
the unital JB$^*$-algebras $M$ and $M(w)$ are not Jordan $^*$-isomorphic. For the concrete example it suffices to consider the transposition on $M_2(\mathbb{C})$, the unit circle $\mathbb{T}$ in $\mathbb{C},$ the JB$^*$-algebra $$M= \{ a:\mathbb{T}\to M_2(\mathbb{C}) \hbox{ continuous} : a(\lambda)^{t} = a(\lambda) \hbox{ for all } \lambda \in \mathbb{T} \}$$ of all continuous functions from $\mathbb{T}$ to the JB$^*$-algebra of all symmetric {\rm(}complex{\rm)} matrices in $M_2(\mathbb{C})$ equipped with the natural Jordan product $a\circ b := \frac12 (ab + b a)$, and the unitary $w : \mathbb{T}\to M_2(\mathbb{C}),$
$w(\lambda) = \left(
 \begin{array}{cc}                                                                                                                                                  \lambda & 0 \\
 0& 1 \\
 \end{array}
\right)$.\smallskip

Let $\mathcal{U}^c (M)$  denote the connected component of $\mathcal{U} (M)$ containing $w$.  We can now conclude that there exists no surjective isometry from $\mathcal{U}^0 (M)$ onto $\mathcal{U}^c (M)$. Otherwise, Theorem \ref{t surj isom secondary components} would imply that $M$ and $M(w)$ are Jordan $^*$-isomorphic, which is impossible.\smallskip

Consider now the unitary $v= w^2\in M$. Let $\mathcal{U}^{\tilde{c}} (M)$ denote the connected component of $\mathcal{U} (M)$ containing $v$. Since $v$ is the square of $w$ and the mapping $U_w: M\to M(v)$ is a Jordan $^*$-isomorphism mapping $\11b{M}$ to $v$ {\rm(}cf. Remark \ref{r unitaries which are squares}{\rm)}, it follows that $U_w \left( \mathcal{U}^0 (M) \right) = \mathcal{U}^{\tilde{c}} (M).$ 
\end{remark}\smallskip

\subsection{Extensibility of surjective isometries between Jordan unitary sets}\label{subsec: extensibility}\ \\

We have already noted the necessity of working with surjective real linear isometries when studying surjective isometries between connected components of the unitary sets of two unital JB$^*$-algebras. More concretely, contrary to the conclusions in the setting of unital C$^*$-algebras, Jordan $^*$-isomorphisms (for the original Jordan product and involution) are not enough to describe these maps (cf. Theorem \ref{t surj isom secondary components} and Remark \ref{r existence of unitaries with non isomorphic nor isometric connected components in the Jordan setting}).\smallskip

As observed by O. Hatori and L. Moln{\'a}r in the setting of commutative unital C$^*$-algebras (see \cite[Corollary 8]{HatMol2014}), and subsequently by Hatori in the more general case of unital C$^*$-algebras $A$ and $B$, a surjective isometry $\Delta: \mathcal{U} (A)\to \mathcal{U} (B)$ need not admit an extension to a surjective real linear isometry from $A$ onto $B$ (cf. \cite[Corollary 5.1]{Hatori14}). We are now in a position to consider the extendibility of a surjective isometry between the unitary sets of two unital JB$^*$-algebras. We include here a new argument based on a variant of the Russo--Dye theorem for unital JB$^*$-algebras due to J.D.M. Wright and M.A. Youngson \cite{WriYou77}.

\begin{corollary}\label{c extendibility of surjective isometries}  Let $\Delta :  \mathcal{U} (M)\to \mathcal{U} (N)$ be a surjective isometry between the unitary sets of two unital JB$^*$-algebras. Let $\{\mathcal{U}^0 (M)\}\cup \{ \mathcal{U}^{j} (M) : j\in \Lambda\}$ be the collection of all connected components of $\mathcal{U} (M)$. For each $j\in \Lambda\cup \{0\}$, let $T_j : M\to N$ denote the surjective real linear isometry extending the mapping $\Delta|_{\mathcal{U}^{j} (M)} : \mathcal{U}^{j} (M)\to \Delta(\mathcal{U}^{j} (M))$ whose existence is assured by Theorem \ref{t surj isom secondary components}. Then the following statements are equivalent:
\begin{enumerate}[$(a)$]\item $\Delta$ admits an extension to a surjective {\rm(}real linear{\rm)} isometry from $M$ onto $N$;
\item The extensions $T_{j_1}$ and $T_{j_2}$ coincide for all $j_1,j_2\in \Lambda\cup \{0\}$.
\end{enumerate}
\end{corollary}

\begin{proof} $(a)\Rightarrow(b)$ Suppose there exists a surjective isometry $T: M\to N$ such that $T|_{\mathcal{U} (M)} = \Delta$. Fix $j\in \Lambda\cup \{0\}$, by hypothesis, $T(u) = \Delta (u) = T_j (u)$ for all $u\in  \mathcal{U}^{j} (M)$.\smallskip

Fix $u_0\in \mathcal{U}^{j} (M)$. Let us consider the unital JB$^*$-algebra given by the $u_0$-isotope $M(u_0)$. By the variant of the Russo--Dye theorem proved in \cite[Corollary]{WriYou77}, the closed unit ball of the space $M = M(u_0)$ is the closed convex hull of the set
$$\{ \exp_{u_0} (i h) : h\in (M(u_0))_{sa} \}.$$ By applying Remark \ref{r connected components which are not the pricipal one} we deduce that the latter set is contained in the principal component of $\mathcal{U} (M(u_0))$ which is precisely $\mathcal{U}^{j} (M)$. It follows from the real linearity of $T$ and $T_j,$ combined with the fact that these two maps coincide on $\mathcal{U}^{j} (M)$, that $T=T_j$ on the closed unit ball of $M$, and hence on the whole $M$.\smallskip

$(b)\Rightarrow(a)$ This implication is easy, it suffices to take $T= T_j$ for any $j\in \Lambda\cup \{0\}$.
\end{proof}

It is not hard to find surjective isometries between sets of unitary elements in two unital JB$^*$-algebras which are not extendable to surjective real linear isometries. Namely, as in Remark \ref{r existence of unitaries with non isomorphic nor isometric connected components in the Jordan setting} we can find a unital JB$^*$-algebra $M$ and a unitary element $w_0$ such that the connected components $\mathcal{U}^0 (M)$ and $\mathcal{U}^c (M)$ are not isometric (where $\mathcal{U}^c (M)$ is the connected component of $\mathcal{U} (M)$ containing $w_0$). Let us pick two surjective linear isometries $T_1, T_2 : M\to M$ with $T_1\neq T_2$ and $T_j (\mathcal{U}^0 (M)) = \mathcal{U}^0 (M)$ for all $j=1,2$ --consider, for example, $T_j= U_{u_j}$ with $u_1\neq u_2$ in $\mathcal{U}^0 (M)$. We define a mapping $\Delta : \mathcal{U} (M) \to \mathcal{U} (M)$ given by $\Delta (u) = T_1(u)$ if $u\in \mathcal{U}^0 (M)$ and $\Delta (u) = T_2 (u)$ otherwise. The surjectivity of $\Delta$ is clear because $T_j (\mathcal{U}^0 (M)) = \mathcal{U}^0 (M)$, and thus $T_2 (\mathcal{U} (M)\backslash \mathcal{U}^0 (M)) = \mathcal{U} (M)\backslash \mathcal{U}^0 (M)$. Since, by Theorems \ref{t principal component as product of exp} and \ref{r connected components which are not the pricipal one}, given $u,v$ belonging to two different connected components of $\mathcal{U} (M)$ we have $\|u-v\| =2$, it follows from the definition of $\Delta$ that it is an isometry. Obviously, $\Delta$ is non-extendable by Corollary \ref{c extendibility of surjective isometries}.\smallskip

\begin{proposition}\label{p T and T0} Let $T,T_0 : M\to N$ be to real linear surjective isometries between two unital JB$^*$-algebras. Let us assume that $T_0$ is a real linear Jordan $^*$-isomorphism, and suppose that $U_{T_0(u_0)} = U_{T(u_0)}$ for all $u_0\in \mathcal{U}^0 (M)$. Then there exists a central projection $p\in M$ such that $T|_{M\circ p\equiv M_2(p)} = T_0|_{M\circ p\equiv M_2(p)}$ and $T|_{M\circ (\textbf{1}-p)\equiv M_2(\textbf{1}-p)} = - T_0|_{M\circ (\textbf{1}-p)\equiv M_2(\textbf{1}-p)}$.
\end{proposition}

\begin{proof} Assume first that $T$ is unital and hence a real linear Jordan $^*$-isomorphism (cf. the comments in page \pageref{eq unital triple hom are Jordan star hom}). It follows from the assumptions that given $h\in M_{sa}$ we have $$U_{T_0(e^{it h})} (\textbf{1}) = U_{T(e^{it h})} (\textbf{1}),$$ for all $t\in \mathbb{R}$. Taking derivatives at $t=0$ in the above equality we arrive at $$ 2 T_0( i h)= 2 U_{T_0( i h),\textbf{1}} (\textbf{1}) = 2 U_{T(i h),\textbf{1}} (\textbf{1})= 2 T(i h),$$ which proves that $T_0( i h) = T( i h)$ for all $h\in M_{sa}$. We should recall that $T$ and $T_0$ are merely real linear maps. Having in mind that $T$ and $T_0$ are real linear Jordan $^*$-isomorphisms we deduce that $$- T(h^2) = T(ih)^2 = T_0(ih)^2 = - T_0(h^2), \hbox{ for all } h\in M_{sa},$$ witnessing that $T$ and $T_0$ coincide on positive elements of $M$. Since every element in $M_{sa}$ writes as the orthogonal sum of two positive elements in $M$, we conclude that $T$ and $T_0$ coincide on $M_{sa}$, and for each $x\in M$ we set $x = h + i k$ with $h,k\in M_{sa}$ to prove that $$T(x) = T(h +i k) = T(h) + T(i k) = T_0(h) + T_0(i k) = T_0(x).$$ We have therefore proved that $T= T_0$ in this case.\smallskip

By hypothesis, $Id= U_{T_0(\textbf{1})} = U_{T(\textbf{1})},$ and hence $T(\textbf{1})^* = U_{T(\textbf{1})}(T(\textbf{1})^*) = T(\textbf{1}),$ witnessing that $T(\textbf{1})$ is a symmetry in $N$, that is, there exist two orthogonal projections $q_1,q_2\in N$ such that $T(\textbf{1}) = q_1 - q_2$ and $q_1 + q_2=\textbf{1}$. Since $U_{q_1} + U_{q_2} - 2 U_{q_1,q_2} =U_{q_1-q_2} = Id =U_{q_1+ q_2} = U_{q_1} + U_{q_2} + 2 U_{q_1,q_2},$ it can be easily deduced that $U_{q_1,q_2} =0$ and hence $N = N_2(q_1) \oplus^{\infty} N_2(q_2)$, which implies that $q_1$ and $q_2$ are central projections.\smallskip

Since $T_0$ is a real linear Jordan $^*$-isomorphism, the elements $p_1 = T_0^{-1} (q_1)$ and $p_2 = T_0^{-1} (q_2)$ are two orthogonal central projections in $M$ with $p_1 + p_2 =\textbf{1}$ and $M = M_2(p_1) \oplus^{\infty} M_2(p_2)$. Under these circumstances, every unitary $u$ in $M$ is of the form $u = u_1 + u_2$ with $u_j \in M_2(p_j)$, and $\mathcal{U}^0 (M)$ identifies with $\mathcal{U}^0 (M_2(p_1))\times \mathcal{U}^0 (M_2(p_2))$. It is clear from the hypotheses that $T_0 $ is the orthogonal sum of the two real linear Jordan $^*$-isomorphisms $T_0|_{M_2(p_1)} : M_2(p_1) \to N_2(q_1)$ and $T_0|_{M_2(p_2)} : M_2(p_2)\to N_2(q_2)$. Moreover, since $T$ is a surjective real linear isometry (and hence a real linear triple isomorphism, cf. \cite[Corollary 3.2]{Da} or \cite[Corollary 3.4]{FerMarPe}) with $T(\textbf{1}) = q_1 - q_2$. Therefore, the mapping $T: M\to N(q_1-q_2)$ is a real linear Jordan $^*$-isomorphism, and $T$ writes as the orthogonal sum of the two real linear Jordan $^*$-isomorphisms $T|_{M_2(p_1)} : M_2(p_1) \to N_2(q_1)$ and $T|_{M_2(p_2)} : M_2(p_2)\to N_2(-q_2)$.\smallskip

We claim that the pairs $T_0|_{M_2(p_1)}, T|_{M_2(p_1)} : M_2(p_1) \to N_2(q_1)$ and $-T_0|_{M_2(p_2)},$ $T|_{M_2(p_2)} : M_2(p_2)\to N_2(-q_2)$ (which are clearly unital) satisfy the hypothesis assumed for the pair $T_0,T$. Namely, given $u_j\in \mathcal{U}^0(M_2(p_j))$, the element $u = u_1 + u_2$ is a unitary in $\mathcal{U}^0 (M),$ and by our hypotheses and the fact that $M = M_2(p_1) \oplus^{\infty} M_2(p_2)$ and $N = N_2(q_1) \oplus^{\infty} N_2(q_2)$ we obtain that $$ U_{T(u_1)} + U_{T(u_2)}=  U_{T(u_1+u_2)} = U_{T_0(u_1+u_2)} = U_{T_0(u_1)} + U_{T_0(u_2)},$$ which implies that $U_{T_0(u_1)} =  U_{T(u_1)}$ and $U_{T_0(u_2)} = U_{T(u_2)}$, because these maps have orthogonal ranges and supports. This finishes the proof of the claim.\smallskip

We can therefore apply the conclusion in the first paragraph of this proof to the pairs $T_0|_{M_2(p_1)}$ and $T|_{M_2(p_1)}$ and $-T_0|_{M_2(p_2)}$ and $T|_{M_2(p_2)}$ to obtain that $T_0|_{M_2(p_1)}=T|_{M_2(p_1)}$ and $-T_0|_{M_2(p_2)}=T|_{M_2(p_2)}$. The proof concludes by taking $p = p_1$.
\end{proof}

We have seen in the previous results that the extendibility of a surjective isometry $\Delta$ between the sets of unitary elements in two unital JB$^*$-algebras $M$ and $N$ is not always an easy task. Actually, the extendibility is, in general, hopeless. It is natural to ask whether an additional algebraic hypothesis can be added to guarantee a linear extension. The instability of the set of unitaries under Jordan products induces us to discard the preservation of the Jordan product as extra algebraic hypothesis. In view of the properties of surjective real linear isometries between JB$^*$-algebras and the conclusion in \cite[Proposition 3.9]{Pe2020}, it seems natural to assume the following extra premise: $$\Delta U_u (v^*) = U_{\Delta(u)} \left(\Delta(v)^* \right), \hbox{ for all } u,v\in \mathcal{U} (M).$$ We have been able to find a positive answer under some extra suppositions on the connected components of the unitary set. We shall only remark that the extra assumption is weaker than the property asserting that every unitary element admits a unitary square root.

\begin{proposition}\label{p extra algebraic hypotheses} Let $\Delta :  \mathcal{U} (M)\to \mathcal{U} (N)$ be a surjective isometry between the unitary sets of two unital JB$^*$-algebras. Let $\{\mathcal{U}^0 (M)\}\cup \{ \mathcal{U}^{j} (M) : j\in \Lambda\}$ be the collection of all connected components of $\mathcal{U} (M)$. Assume that for each $j\in \Lambda$ there exists $u_j \in \mathcal{U} (M)$ with $u_j^2 \in \mathcal{U}^{j} (M)$. Suppose additionally that \begin{equation}\label{eq preservation of quadratic products extra hypo} \Delta U_u (v^*) = U_{\Delta(u)} \left(\Delta(v)^* \right), \hbox{ for all } u,v\in \mathcal{U} (M).
\end{equation} Then $\Delta$ admits an extension to a surjective {\rm(}real linear{\rm)} isometry from $M$ onto $N$.
\end{proposition}

\begin{proof} Up to replacing the original Jordan product and involution on $N$ by the one in $\Delta(\textbf{1})$-isotope, we can assume that $\Delta$ is unital, and hence $\Delta (\mathcal{U}^0 (M)) = \mathcal{U}^0 (N)$.  Let $T_j : M\to N$ denote the surjective real linear isometry extending the mapping $\Delta|_{\mathcal{U}^{j} (M)} : \mathcal{U}^{j} (M)\to \Delta(\mathcal{U}^{j} (M))$ whose existence is assured by Theorem \ref{t surj isom secondary components}.\smallskip

Fix $j,k\in \Lambda$ satisfying $\left(\mathcal{U}^j (M)\right)^* = \mathcal{U}^k(M)$. For each $u_j\in \mathcal{U}^j (M)$, it follows from \eqref{eq preservation of quadratic products extra hypo} that $$T_k (u_j^*) = \Delta (u_j^*) = \Delta U_{\textbf{1}} (u_j^*) =  U_{\textbf{1}} (\Delta(u_j)^*) = T_j (u_j)^*,$$ which proves that $T_k (u_j^*) = T_j (u_j)^*$ for all $u_j \in \mathcal{U}^j (M)$. We fix $w_j\in \mathcal{U}^j (M)$. Since, by the Wright-Youngson-Russo-Dye theorem \cite[Corollary]{WriYou77}, the convex hull of the set $\{ \exp_{w_j} (i h) : h\in (M(w_j))_{sa} \}(\subset \mathcal{U}^j (M) = \mathcal{U}^0 (M(w_j)))$ is norm dense in the closed unit ball of $M(w_j)$, it follows from the linearity and continuity of $T_k$ and $T_j$ that \begin{equation}\label{eq Tj and Tk coincide} T_k (x^*) = T_j (x)^* \hbox{ for all } x \in M.
\end{equation}

Take now an arbitrary $u_0 \in \mathcal{U}^0(A)$. By \eqref{eq inner main other main} we have $U_{u_0} (u_j^*)\in \mathcal{U}^k(M),$ and hence by \eqref{eq preservation of quadratic products extra hypo} we get $$\begin{aligned}
U_{T_k(u_0)} (T_j(u_j)^*) &= U_{T_k(u_0)} (T_k(u_j^*)) = T_k U_{u_0} (u_j^*) = \Delta U_{u_0} (u_j^*) = U_{\Delta(u_0)} \left(\Delta(u_j^*)\right) \\
&= U_{T_0(u_0)} \left( T_k(u_j^*)\right)= \hbox{\eqref{eq Tj and Tk coincide}}= U_{T_0(u_0)} \left( T_j(u_j)^*\right)
\end{aligned},$$ for all $u_j\in \mathcal{U}^j (M)$. A new application of the Wright-Youngson-Russo-Dye theorem as in the previous paragraph --together with the surjectivity of $T_j$-- implies that $U_{T_k(u_0)} = U_{T_0(u_0)}$ for all $u_0 \in \mathcal{U}^0(M)$. Proposition \ref{p T and T0} assures the existence of a central projection $p_k\in M$ such that \begin{equation}\label{eq existence of pk} T_k = T_0|_{M_2(p_k)} \oplus T_0|_{M_2(\textbf{1}- p_k)}.
\end{equation}

Fix an arbitrary $l\in \Lambda$. By assumptions, there exists $u_l\in \mathcal{U} (M)$ with $u_l^2 \in \mathcal{U}^{l} (M)$. We can assume that $u_l\in \mathcal{U}^{j} (M)$ and $u_l^*\in \mathcal{U}^{k} (M)$ for some $j,k\in \Lambda$. By hypotheses, \begin{equation}\label{eq Tl and Tj ul square} T_l (u_l^2) = \Delta (u_l^2) = \Delta U_{u_l} (\textbf{1}) =  U_{\Delta(u_l)} \Delta (\textbf{1}) = \Delta(u_l)^2 = T_j (u_l)^2.
 \end{equation}By considering the central projections $p_l$ and $p_j$ given by \eqref{eq existence of pk}, by orthogonality, we have $$ T_j (u_l)^2 = \left( T_0 (u_l\circ p_j) - T_0(u_l\circ (\textbf{1}-p_j))\right)^2 = T_0 (u_l^2 \circ p_j) + T_0(u_l^2\circ (\textbf{1}-p_j) ) = T_0 (u_l^2)$$ and $$ T_l (u_l^2) = T_0 (u_l^2\circ p_l) - T_0(u_l^2\circ (\textbf{1}-p_l)),$$ which combined with \eqref{eq Tl and Tj ul square} give $T_0(u_l^2\circ (\textbf{1}-p_l))=0,$ and hence $u_l^2\circ (\textbf{1}-p_l)=0$ (by the injectivity of $T_0$). Having in mind that $p_l$ is a central projection, we have $$0 = u_l^2\circ (\textbf{1}-p_l) = U_{u_l} (\textbf{1}-p_l),$$ witnessing that $\textbf{1}-p_l=0,$ and consequently, $T_l = T_0$. The arbitrariness of $l\in \Lambda$ and Corollary \ref{c extendibility of surjective isometries} give the desired statement.
\end{proof}

The result in the previous proposition is new even in the setting of unital C$^*$-algebras.

\medskip\medskip

\noindent\textbf{Acknowledgements} We would like to thank Prof. Lajos Moln{\'a}r for encouraging us to explore this problem. \smallskip

\noindent First and fifth authors partially supported by the Spanish Ministry of Science, Innovation and Universities (MICINN) and European Regional Development Fund project no. PGC2018-093332-B-I00, Programa Operativo FEDER 2014-2020 and Consejer{\'i}a de Econom{\'i}a y Conocimiento de la Junta de Andaluc{\'i}a grant numbers A-FQM-242-UGR18 and FQM375. First author partially supported by EPSRC (UK) project ``Jordan Algebras, Finsler Geometry and Dynamics'' ref. no. EP/R044228/1. Second author partially supported by JSPS KAKENHI Grant Number JP 21J21512. 
Fourth author partially supported by JSPS KAKENHI (Japan) Grant Number JP 20K03650.

\end{document}